\documentclass[a4paper]{amsart}
\usepackage{multicol,ifthen}
\begin{document}

\title{Bi- and trilinear Schr\"odinger estimates in one space dimension with applications
to cubic NLS and DNLS}

\author{Axel Gr\"unrock \\Fachbereich C:  Mathematik/Naturwissenschaften\\ Bergische
Universit\"at Wuppertal \\ Gau{\ss}stra{\ss}e 20 \\ D-42097 Wuppertal \\ Germany \\
e-mail Axel.Gruenrock@math.uni-wuppertal.de}

\date{}

\newcommand{\F}{\mathcal{F}}
\newcommand{\R}{\mathbb R}
\newcommand{\N}{\mathbb N}
\newcommand{\C}{\mathbb C}  
\newcommand{\h}[2]{\mbox{$ \widehat{H^{#1}_{#2}}$}}
\newcommand{\hh}[3]{\mbox{$ \widehat{H^{#1}_{#2, #3}}$}} 
\newcommand{\n}[2]{\mbox{$ \| #1\| _{ #2} $}} 
\newcommand{\x}{\mbox{$X^r_{s,b}$}} 
\newcommand{\xx}{\mbox{$X_{s,b}$}}
\newcommand{\X}[3]{\mbox{$X^{#1}_{#2,#3}$}} 
\newcommand{\XX}[2]{\mbox{$X_{#1,#2}$}}
\newcommand{\q}[2]{\mbox{$ {\| #1 \|}^2_{#2} $}}
\newcommand{\e}{\varepsilon}
\pagestyle{plain} 

\newtheorem{lemma}{Lemma} 
\newtheorem{kor}{Corollary} 
\newtheorem{satz}{Theorem}
\newtheorem{prop}{Proposition}

\begin{abstract}
The Fourier transforms of the products of two respectively three solutions of the free Schr\"odinger equation in one space dimension are
estimated in mixed and, in the first case, weighted $L^p$ - norms. Inserted into an appropriate variant of the Fourier restriction norm
method, these estimates serve to prove local well-posedness of the Cauchy problem for the cubic nonlinear Schr\"odinger (NLS) equation with
data $u_0$ in the function spaces $\widehat{L^r_x}:=\widehat{H^r_0}$, where for $s \in \R$ the spaces $\widehat{H^r_s}$ are defined by their
norms
\[\n{u_0}{\widehat{H^r_s}}:=\n{\langle \xi \rangle^s \widehat{u_0}}{L^{r'}_{\xi}}, \,\,\,\,\,\frac{1}{r}+\frac{1}{r'}=1.\]
Similar agruments, combined with a gauge transform, lead to local well-posedness of the Cauchy problem for the derivative nonlinear Schr\"odinger
(DNLS) equation with data $u_0 \in \widehat{H^r_{\frac{1}{2}}}$. In the local result on cubic NLS the parameter $r$ is allowed in the full subcritical
range $1 < r < \infty$, while for DNLS we assume $1<r\le 2$. In the special case $r=2$ both results coincide with the optimal ones on the $H^s$ - scale. 

Furthermore, concerning the cubic NLS equation, it is shown by a decomposition argument that the local solutions extend globally, provided 
$2 \ge r > \frac{5}{3}$.
\end{abstract}

\maketitle

\section{Introduction and main results}

The Cauchy problem for the cubic nonlinear Schr\"odinger equation on the line
\begin{equation}\label{101}
iu_t + u_{xx} + |u|^2 u = 0, \hspace{2cm} u(0,x)=u_0(x) ,\,\,\,x\in \R
\end{equation}
is known to be globally well-posed for data in the classical Sobolev spaces $H^s$, if $s \ge 0$, and locally ill-posed in the sense that the
mapping data upon solution fails to be uniformly continuous, if $s<0$. The well-posedness result goes back to Y. Tsutsumi \cite{T87} (see also
Cazenave and Weissler \cite{CW90}), while ill-posedness below $L^2_x$ has been shown by Kenig, Ponce and Vega in \cite{KPV01}. The ``criticality''
of $L^2_x$ (in the sense that well-posedness holds for $H^s$ - data, iff $H^s \subset L^2_x$, i. e. iff $s \ge 0$) can be explained heuristically
by Galilean invariance, see the introduction of \cite{KPV01}. On the other hand scaling considerations suggest local well-posedness for a larger
class of data\footnote{From the scaling point of view the critical Sobolev index for the IVP (\ref{101}) is $s=-\frac{1}{2}$.}, and in fact local and
even global well-posedness of (\ref{101}) for data with an infinite $L^2_x$ - norm has been demonstrated by Vargas and Vega in \cite{VV01}. 
Inspired by these results as well as by the work of Cazenave, Vega and Vilela \cite{CVV01} we consider the Cauchy problem (\ref{101}) with data
$u_0$ in the space $\h{r}{s}$, which is defined by the norm
\[\n{u_0}{\h{r}{s}} := \n{\langle \xi \rangle ^s\widehat{u_0}}{L^{r'}_{\xi}},\hspace{1cm}\frac{1}{r}+\frac{1}{r'}=1,\]
where $\langle \xi \rangle ^s = (1+|\xi|^2)^{\frac{s}{2}}$. For $s=0$ we will write $\widehat{L^r_x}$ instead of $\h{r}{0}$. Concerning the Cauchy problem (\ref{101}) we will show local well-posedness in
$\h{r}{s}$ for $s\ge 0$ and $1<r<\infty$ (see Theorem \ref{t1} below). Observe that $\h{2}{s}= H^s$, so for $r=2$ this coincides with the optimal local $H^s$ - result.
Furthermore, $\h{r}{s}$ scales like $H^{\sigma}$, if $s-\frac{1}{r}=\sigma-\frac{1}{2}$, hence from the scaling point of view we obtain an improvement
by pushing down $r$ from $2$ to $1+$, where - for $s=0$ - we almost reach the scaling line $s-\frac{1}{r}=-1$. In this setting, the case $(s,r)=(0,1)$ 
becomes critical (with respect to both, the Galilean and the scaling transformations) and must be left as an open problem.\\

As long as $\frac{4}{3} < r \le 2$, our result can be obtained quite easily by using the linear estimate
\begin{equation}\label{102}
\n{e^{it\partial_x ^2} u_0}{L^{3r}_{xt}} \le c \n{{u_0}}{\widehat{L^r_x}},
\end{equation}
($\frac{4}{3} < r \le \infty$), which goes back to Fefferman and Stein \cite{F70}. This is already contained in the arguments of \cite{VV01} and
\cite{CVV01}, see also \cite[Proposition 1.1.]{G04}. Unfortunately, the estimate (\ref{102}) fails for $r \le \frac{4}{3}$. To overcome this difficulty,
we use bi- and trilinear estimates for free solutions of the linear Schr\"odinger equation. More precisely, for $u=e^{it\partial_x ^2}u_0$, $v=e^{it\partial_x ^2}v_0$
and $w=e^{-it\partial_x ^2}w_0$ we estimate $I^{\frac{1}{p}}(vw)$ ($I$ being the Riesz potential operator of order $-1$) and the product $uvw$ in
the mixed space-time norms
\[\n{f}{\widehat{L_x^q}(\widehat{L_t^p})}:= \left( \int \Big{(} \int |\widehat{f}(\xi, \tau)|^{p'} d \tau \Big{)}^{\frac{q'}{p'}}d\xi \right)^{\frac{1}{q'}}, \,\,\,\,\,\frac{1}{q}+\frac{1}{q'}=\frac{1}{p}+\frac{1}{p'}=1.\]
(Here $f \in \mathcal{S}' (\R^2)$ depends on the space variable $x$ and the time variable $t$, $\widehat{f}$ is its Fourier transform with respect to space
and time and $(\xi, \tau)$ denote the variables conjugate to $(x,t)$. When $p=q$, we will write for short $\widehat{L_{xt}^p}$ instead of $\widehat{L_x^p}(\widehat{L_t^p})$.)\\

These multilinear estimates, whose precise statement and proof is content of section 2, are then inserted into the framework of Bourgain's 
Fourier restriction norm method (see \cite{B93}) respectively into its generalization to non $L^2$
-based spaces developed by the author in \cite{G04}. We recall the function spaces
\[\x:= \{f \in \mathcal{S'}(\R ^{2}): \n{f}{\x}< \infty\},\]
where $s,b \in \R$, $1 \le r \le \infty$, $\frac{1}{r} + \frac{1}{r'}=1$ and
\[\n{f}{\x}:= \left(\int d \xi d \tau \langle \xi \rangle^{sr'}\langle \tau + \xi ^2\rangle^{br'} |\hat{f}(\xi , \tau)|^{r'} \right) ^{\frac{1}{r'}}\]
with the usual modification for $r=1$, as well as the time restricted spaces
\[\x(\delta) := \{f = \tilde{f}|_{[-\delta,\delta] \times \R} : \tilde{f} \in \x\}\]
endowed with the norm
\[\n{f}{\x(\delta)}:= \inf \{ \n{\tilde{f}}{\x} : \tilde{f}|_{[-\delta,\delta] \times \R^n} =f\}  .\]
For $r=2$ these are the fuction spaces $\xx$ (respectively $\xx(\delta)$) introduced by Bourgain in \cite{B93} in the study of initial value problems. So in this case we shall omit the index $r$.\\

In this framework, concerning the Cauchy problem
\begin{equation}\label{103}
iu_t + u_{xx} = N(u), \hspace{2cm}u(0)=u_0 \in \h{r}{s},
\end{equation}
with a general nonlinearity $N$ depending on $u$ and its derivatives, the following local well-posedness theorem holds true:
\setcounter{satz}{-1}
\begin{satz}\label{t0} Assume that for given $s \in \R$, $r \in (1, \infty)$ there exist $b > \frac{1}{r}$ and $b' \in (b-1,0]$, such that the estimates
\begin{equation}\label{104}
\n{N(u)}{\X{r}{s}{b'}} \le C( \|u\|_{\x})\|u\|_{\x}
\end{equation}
and
\begin{equation}\label{105}
\n{N(u) - N(v)}{\X{r}{s}{b'}} \le C (\|u\|_{\x} + \|v\|_{\x})\|u - v\|_{\x}
\end{equation}
are valid with a continuous and nondecreasing Function $C:\R_0^+ \rightarrow \R_0^+$. Then there exist $\delta = \delta (\n{u_0}{\h{r}{s}}) > 0$ and a unique solution $u \in \x(\delta)$ of (\ref{103}). This solution is persistent and the mapping S: $u_0 \mapsto u$, $\h{r}{s} \rightarrow \x (\delta _0)$ (data upon solution) is locally Lipschitz continuous for any $\delta _0 \in (0,\delta )$.
\end{satz}
See \cite[Theorem 2.3.]{G04}. The replacement of the special function $C(t)=ct^{\alpha -1}$ in that theorem by an arbitrary continuous and nondecreasing
function is obvious. It should be remarked here, that by a solution of (\ref{103}) we always mean a solution of the corresponding integral equation
\begin{equation}\label{106}
u(t)= e^{it\partial_x^2}u_0 - i \int_0^t e^{i(t-s)\partial_x^2}N(u)(s)ds.
\end{equation}
Moreover, let us for further reference recall the two linear estimates needed in the proof of Theorem \ref{t0}, these are
\begin{equation}\label{107}
\n{e^{it\partial_x^2}u_0}{\x(\delta)} \le c \n{u_0}{\h{r}{s}}
\end{equation}
and, provided $1<r<\infty$ and $b'+1 \ge b \ge 0 \ge b' > -\frac{1}{r'}$,
\begin{equation}\label{108}
\n{\int_0^t e^{i(t-s)\partial_x^2}F(s)ds}{\x(\delta)} \le c \delta^{1-b+b'} \n{F}{\X{r}{s}{b'}(\delta)},
\end{equation}
see (2.17) and Lemma 2.2. in \cite{G04}. After these preparations our result concerning (\ref{101}) simply reads:
\begin{satz}\label{t1} 
Let $s\ge 0$ and $1<r<\infty$. Then the Cauchy problem (\ref{101}) is locally well-posed in $\h{r}{s}$ in the sense of Theorem \ref{t0}.
\end{satz}

When combined with the gauge transform
\[Gf(x):=e^{-i\int_{- \infty}^x |f(y)|^2 dy}f(x)\]
our arguments also apply to the Cauchy problem for the derivative nonlinear Schr\"odinger equation (DNLS) in one space dimension
\begin{equation}\label{109}
iu_t + u_{xx} = i (|u|^2 u)_x, \hspace{1cm} u(x,0)=u_0(x), \,\,\,\,\,\,x\in\R.
\end{equation}
This problem has been shown to be locally well posed for $H^s$-data, $s \ge \frac{1}{2}$, by Takaoka in 1999 \cite{T99},
where he improved earlier results of Hayashi and Ozawa \cite{HO94}. On the $H^s$-scale, the $H^{\frac{1}{2}}$-
result is optimal; in fact, ill-posedness in the $C^0$-uniform sense has been demonstrated by Biagioni and 
Linares in 2001 \cite{BL01} using an appropriate counterexample. On the other hand the critical scaling exponent for equation (\ref{109})
is $s=0$. Again there is a gap of half a derivative between the optimal local well-posedness result in $H^s$ and the scaling prediction.
Proving local well-posedness of (\ref{109}) in $\h{r}{s}$ for $s \ge \frac{1}{2}$, $2 \ge r > 1$, we can close this gap at least partially.
In order to do so, we follow Takaoka and consider first the gauge equivalent problem
\begin{equation}\label{110}
iv_t + v_{xx} + i v^2 \overline{v}_x + \frac{1}{2}|v|^4v=0, \hspace{1cm} v(x,0)=v_0(x), \,\,\,\,\,\,x\in\R,
\end{equation}
for which we can show
\begin{satz}\label{t2} 
Let $s\ge \frac{1}{2}$ and $1<r \le 2$. Then the Cauchy problem (\ref{110}) is locally well-posed in $\h{r}{s}$ in the sense of Theorem \ref{t0}.
\end{satz}

The nonlinear estimates necessary for Theorem \ref{t2} (cf. (\ref{104}) and (\ref{105})) are proven in section 3.1, see Lemma \ref{l30} and Lemma \ref{l31} below.
In a second step involving the gauge transform the following result concerning (\ref{109}) is obtained:
\begin{satz}\label{t3}
Let $s\ge \frac{1}{2}$ and $1<r \le 2$. Then for $u_0 \in \h{r}{s}$ there exist $b>\frac{1}{r}$, $\delta = \delta (\n{u_0}{\h{r}{s}}) > 0$ and a unique solution
\[u \in G^{-1}(\x(\delta)) \subset C([-\delta,\delta], \h{r}{s})\]
of (\ref{109}). For any $\delta_0 \in (0,\delta)$ the mapping $u_0 \mapsto u$, $\h{r}{s} \rightarrow C([-\delta_0,\delta_0], \h{r}{s})$ (data 
upon solution) is locally Lipschitz continuous.
\end{satz}

(Strictly speaking, the uniqueness statement in Theorem \ref{t3} is to be understood in the following sense: If $u' \in G^{-1}(\x(\delta))$ is a further solution
of (\ref{109}) such that $Gu$ solves (\ref{110}) with $v_0=Gu_0$, then $u'=u$. The additional hypothesis - ``such that $Gu$ solves (\ref{110})'' -
seems to be somewhat weak and artificial, but I cannot see, how to remove it unless by changing the notion of a solution. For example, this hypothesis
is always satisfied, if there is a sequence $(u_n)_{n \in \N}$ of smooth solutions of (\ref{109}) with $Gu_n \rightarrow G_u$ in $\x(\delta)$.)\\

Theorem \ref{t3} is proved in section 3.2, the main problem here is the continuity of the gauge transform $G : \h{r}{s} \rightarrow \h{r}{s}$,
see Lemma \ref{l32} below. At this point the necessity of the additional hypothesis $r \le 2$ becomes obvious, since the definition of $G$ demands for $\h{r}{s} \subset L^2_x$.\\

In both local results the critical cases - i. e. $(s,r)=(0,1)$ in Theorem \ref{t1} respectively $(s,r)=(\frac{1}{2},1)$ in Theorems \ref{t2} and
\ref{t3} - remain open. Nevertheless, these results are sharp in the sense that for given $r>1$ the initial value problems (\ref{101}) (respectively
(\ref{109})) are ill-posed in $\h{r}{s}$, if $s<0$ (respectively if $s<\frac{1}{2}$). To see this we use the counterexamples from \cite{KPV01} 
concerning the cubic NLS equation and from \cite{BL01} concerning the DNLS equation, see section 4, where the details are discussed.\\

Finally, coming back to the IVP (\ref{101}) in section 5, we show global well-posedness of this problem for data in $\widehat{L^r_x}$ in the
parameter range $2 \ge r > \frac{5}{3}$. More precisely:
\begin{satz}\label{t4}
Let $2 \ge r > \frac{5}{3}$ and $u_0 \in \widehat{L^r_x}$. Then the local solution $u$ of (\ref{101}) obtained in Theorem \ref{t1} extends globally.
Moreover, the difference
\[z(t):=u(t)-e^{it\partial^2_x}u_0\]
belongs to $L_x^2$ for all $t>0$ and satisfies the estimate
\[\n{z(t)}{L_x^2} \le c \langle t \rangle ^{\frac{r'-2}{10-4r'}+}.\]
\end{satz}
(Here and below we write $x\pm$ to denote $x\pm\e$ for arbitrarily small $\e>0$.) To prove this theorem we use a data-decomposition argument as
introduced in \cite[\S 7]{B98}. In connection with the Cauchy problem (\ref{101}) this type of argument has already been used by Vargas and Vega, see
section 2 of \cite{VV01}. The only particularity here is the way of splitting the data adapted to the spaces $\widehat{L^r_x}$, which is horizontally
in Fourier space (instead of vertically as in \cite{B98}), see section 5.1.\\

{\bf{Acknowledgement:}} The author wants to thank Luis Vega for his helpful encouragement to deal with these questions, especially to prove a
global result concerning (\ref{101}). He is also grateful to Carlos Kenig, Herbert Koch and Daniel Tataru for inviting him to a workshop at
Oberwolfach last fall, where parts of this material could be presented.
\section{Key estimates}

\begin{lemma}\label{l20} Let $1 \le q \le r_{1,2} \le p \le \infty$ and $\frac{1}{p}+\frac{1}{q}=\frac{1}{r_1}+\frac{1}{r_2}$.
Then, for $u=e^{it\partial_x ^2}u_0$ and $v=e^{-it\partial_x ^2}v_0$, the estimate
\[\n{I^{\frac{1}{p}}(uv)}{\widehat{L^q_x}(\widehat{L^p_t})} \le c \n{u_0}{\widehat{L^{r_1}_x}}\n{v_0}{\widehat{L^{r_{2}}_x}}\]
holds true.
\end{lemma}

Proof: Computing the Fourier-transform first in space and then in time we obtain
\[\mathcal{F}_x I^{\frac{1}{p}}(uv)(\xi,t)= c|\xi|^{\frac{1}{p}} \int_* e^{it(\xi_1^2-\xi_2^2)}\widehat{u_0}(\xi_1)\widehat{v_0}(\xi_2)d\xi_1\]
and
\begin{eqnarray*}
\mathcal{F} I^{\frac{1}{p}}(uv)(\xi,\tau) & = & c|\xi|^{\frac{1}{p}} \int_* \delta(\tau-\xi_1^2+\xi_2^2)\widehat{u_0}(\xi_1)\widehat{v_0}(\xi_2)d\xi_1 \\ 
& = & c|\xi|^{-\frac{1}{p'}} \widehat{u_0}(\frac{\xi}{2}+\frac{\tau}{2 \xi})\widehat{v_0}(\frac{\xi}{2}-\frac{\tau}{2 \xi}),
\end{eqnarray*}
respectively. (Here $\int_*$ is shorthand for $\int_{\xi_1+\xi_2=\xi}$). Hence
\begin{eqnarray*}
\|\mathcal{F} I^{\frac{1}{p}}(uv)(\xi,\cdot)\|^{p'}_{L^{p'}_{\tau}} &=& c|\xi|^{-1} \int d \tau |\widehat{u_0}(\frac{\xi}{2}+\frac{\tau}{2 \xi})\widehat{v_0}(\frac{\xi}{2}-\frac{\tau}{2 \xi})|^{p'} \\
&=& c \int dx |\widehat{u_0}(\frac{\xi + x}{2}) \widehat{v_0}(\frac{\xi - x}{2})|^{p'}\hspace{1,8cm}(x=\frac{\tau}{\xi})\\
&=& c \int dy |\widehat{u_0}(y)\widehat{v_0}(\xi - y)|^{p'}\hspace{2cm}(y=\frac{x+\xi}{2})\\
&=& c |\widehat{u_0}|^{p'}*|\widehat{v_0}|^{p'}(\xi).
\end{eqnarray*}
Now we choose $r'=\frac{q'}{p'}$ ($\ge 1$, since $p \ge q$) and $\rho_{1,2}$ with $\rho'_{1,2}= \frac{r'_{1,2}}{p'}$. Then 
$\frac{1}{r}=\frac{1}{\rho_1}+\frac{1}{\rho_2}$ and, using Young's inequality in the third step, we get
\begin{eqnarray*}
\|\mathcal{F} I^{\frac{1}{p}}(uv)\|_{L^{q'}_{\xi}(L^{p'}_{\tau})} &=& c \left( \int d\xi (|\widehat{u_0}|^{p'}*|\widehat{v_0}|^{p'}(\xi))^{\frac{q'}{p'}}\right)^{\frac{1}{q'}}\\
&=& c \||\widehat{u_0}|^{p'}*|\widehat{v_0}|^{p'}\|^{\frac{1}{p'}}_{L^{r'}_{\xi}}\\
&\le & c \left(\n{|\widehat{u_0}|^{p'}}{L^{\rho'_1}_{\xi}}\n{|\widehat{v_0}|^{p'}}{L^{\rho'_2}_{\xi}} \right)^{\frac{1}{p'}} = c \n{u_0}{\widehat{L^{r_1}_x}}\n{v_0}{\widehat{L^{r_{2}}_x}}
\end{eqnarray*}
$\hfill \Box$

\vspace{0.5cm}

{\bf{Remark:}} As the proof shows, the inequality in Lemma \ref{l20} becomes an \emph{equality}, if $p=q$.

\vspace{0.5cm}

Arguing similarly as in the proof of Lemma 2.1 in \cite{G04} we obtain:

\begin{kor}\label{k20} For $p$, $q$, $r_{1,2}$ as in the previous Lemma and $b_i > \frac{1}{r_i}$ the estimate
\[\n{I^{\frac{1}{p}}(u \overline{v})}{\widehat{L^q_x}(\widehat{L^p_t})} \le c \n{u}{\X{r_1}{0}{b_1}}\n{v}{\X{r_2}{0}{b_2}}\]
is valid.
\end{kor}

{\bf{Remark:}} The case of the above Corollary, where all the H\"older exponents $p,q, r_i$ are equal to $2$ was shown by Bekiranov, Ogawa and Ponce, see \cite[Lemma 3.2]{BOP98}.

\begin{lemma}\label{l21} Let $q>1$, $0 < \frac{1}{r'_{1,2}}<\frac{1}{p'}< \min({\frac{1}{r_0},\frac{1}{r'_1}+\frac{1}{r'_2}})$ and 
$\sum_{i=0}^2 \frac{1}{r_i}=\frac{1}{q}+\frac{2}{p}$. Then for $u=e^{it\partial_x ^2}u_0$, $v=e^{it\partial_x ^2}v_0$ and $w=e^{-it\partial_x ^2}w_0$
we have
\[\n{uvw}{\widehat{L^q_x}(\widehat{L^p_t})} \le c \n{u_0}{\widehat{L^{r_0}_x}}\n{v_0}{\widehat{L^{r_{1}}_x}}\n{w_0}{\widehat{L^{r_{2}}_x}}.\]
\end{lemma}

Proof: The Fourier-transform of the product in the space variable only is
\[\mathcal{F}_x (uvw) (\xi,t)= c \int_* e^{it(\xi_1^2+\xi_2^2-\xi_3^2)}\widehat{u_0}(\xi_1)\widehat{v_0}(\xi_2)\widehat{w_0}(\xi_3)d\xi_1 d \xi_2 ,\]
where now $\int_* =\int_{\xi_1+\xi_2+\xi_3=\xi}$. From this we get for the Fourier-transform in both variables
\[\mathcal{F}(uvw) (\xi,\tau)= c \int_* \delta(\tau - \xi_1^2-\xi_2^2+\xi_3^2)\widehat{u_0}(\xi_1)\widehat{v_0}(\xi_2)\widehat{w_0}(\xi_3)d\xi_1 d \xi_2 .\]
Now the argument $g(\xi_2) :=\tau - \xi_1^2-\xi_2^2+\xi_3^2$ of $\delta$ vanishes, iff $\xi_2 =\frac{\tau + \xi^2- 2\xi \xi_1}{2(\xi - \xi_1)}=:x$, and we have
$|g'(\xi_2)|=2|\xi - \xi_1|$. This gives
\begin{eqnarray}\label{200}
\mathcal{F}(uvw) (\xi,\tau) = c \int \frac{1}{|\xi - \xi_1|}\widehat{u_0}(\xi_1)\widehat{v_0}(x)\widehat{w_0}(\xi-\xi_1-x)d\xi_1 \hspace{1.5cm}\\
 \le  c\left(\int \frac{|\widehat{u_0}(\xi_1)|^p}{|\xi - \xi_1|^{(1-\theta) p}}d\xi_1 \right)^{\frac{1}{p}}
\left(\int |\widehat{v_0}(x)\widehat{w_0}(\xi-\xi_1-x)|^{p'}|\xi - \xi_1|^{-\theta p'}d\xi_1 \right)^{\frac{1}{p'}}, \nonumber
\end{eqnarray}
where $\theta = \frac{3}{p'}-\frac{1}{r'_1}-\frac{1}{r'_2} \in (0,1)$. Taking the $L_{\tau}^{p'}$-norm of both sides, we obtain
\begin{eqnarray*}
\|\mathcal{F}(uvw) (\xi,\cdot)\|_{L_{\tau}^{p'}} \le c \left( |\widehat{u_0}|^p * |\xi|^{(\theta - 1) p}\right)^{\frac{1}{p}} \times \\
\left(\int |\widehat{v_0}(x)\widehat{w_0}(\xi-\xi_1-x)|^{p'}|\xi - \xi_1|^{-\theta p'}d\xi_1 d\tau \right)^{\frac{1}{p'}}.
\end{eqnarray*}
Changing variables ($x$ as above and $y=\xi - \xi_1 -x$) we see that the second factor is equal to
\[c\left( \int |\widehat{v_0}(x)\widehat{w_0}(y)|^{p'}|x+y|^{1-\theta p'}dxdy\right)^{\frac{1}{p'}} \le c  \n{v_0}{\widehat{L^{r_{1}}_x}}\n{w_0}{\widehat{L^{r_{2}}_x}}\]
by the Hardy-Littlewood-Sobolev inequality, requiring $\theta$ to be chosen as above and $1< \theta p' < 2$; $1 < \frac{r'_i}{p'} < \infty$, $i=1,2$;
which follows from our assumptions. It remains to estimate the $L^{q'}_{\xi}$-norm of the first factor, that is
\begin{eqnarray*}
&& \||\widehat{u_0}|^p* |\xi|^{(\theta -1)p}\|^{\frac{1}{p}}_{L_{\xi}^{\frac{q'}{p}}}\\
& \le & c( \n{|\widehat{u_0}|^p}{L_{\xi}^{\frac{r'_0}{p}}} \n{|\xi|^{(\theta -1)p}}{L_{\xi}^{\frac{1}{(1-\theta)p}, \infty}} )^{\frac{1}{p}}\\
& \le & c \n{u_0}{\widehat{L^{r_0}_x}},
\end{eqnarray*}
where the HLS inequality was used again. For its application we need
\[0<(1-\theta)p<1;\,\,\,\,1<\frac{r'_0}{p}< \frac{1}{1-(1-\theta)p}\,\,\,\,\mbox{and}\,\,\,\,\,\frac{p}{q'}=(1-\theta)p-1+\frac{p}{r'_0},\]
which again follows by the assumptions, as can be easily checked.
$\hfill \Box$

\begin{kor}\label{k21} Let $p,q>1$. Assume $p'>q$ or $p=q$. Then, for $u,v$ and $w$ as in Lemma \ref{l21} the estimate
\[\n{uvw}{\widehat{L^q_x}(\widehat{L^p_t})} \le c \n{u_0}{\widehat{L^{q}_x}}\n{v_0}{\widehat{L^{p}_x}}\n{w_0}{\widehat{L^{p}_x}}\]
holds true.
\end{kor}

Proof: We consider the case $p'>q$ first: For $\rho$, $\rho_0$ with $\frac{4}{3}<\rho_0<2<\rho$ and $\frac{1}{\rho_0}+\frac{2}{\rho}=\frac{3}{2}$ we have by H\"older and (\ref{102})
\[\n{uvw}{L^2_{xt}} \le c \n{u_0}{\widehat{L^{\rho_0}_x}}\n{v_0}{\widehat{L^{\rho}_x}}\n{w_0}{\widehat{L^{\rho}_x}}.\]
We define - for $\theta \in (0,1)$ - the H\"older exponents $q_{\theta}$, $p_{\theta}$, $r_{\theta}$ and $r_{0, \theta}$ by the interpolation
conditions
\[\frac{1}{q}=\frac{1-\theta}{q_{\theta}}+\frac{\theta}{2}=\frac{1-\theta}{r_{0, \theta}}+\frac{\theta}{\rho_0}\,\,\,;
\,\,\,\frac{1}{p}=\frac{1-\theta}{p_{\theta}}+\frac{\theta}{2}=\frac{1-\theta}{r_{\theta}}+\frac{\theta}{\rho}.\]
Then, by multilinear interpolation, it is sufficient to show that - for $\theta$ small enough - the exponents $q_{\theta}$ ($p_{\theta}$)
instead of $q$ ($p$) and $r_0 = r_{0, \theta}$, $r_1=r_2=r_{\theta}$ fulfill the assumptions of Lemma \ref{l21}: The identity
\[\frac{1}{q_{\theta}} + \frac{2}{p_{\theta}} = \frac{1}{r_{0,\theta}} + \frac{2}{r_{\theta}}\]
is easily checked. The condition $q_{\theta} >1$ becomes $\frac{1}{q}<1-\frac{\theta}{2}$ (i). We have $0<\frac{1}{r'_{\theta}}$, iff
$\frac{1}{p}<1-\frac{\theta}{\rho '}$ (ii) and $\frac{1}{r'_{\theta}}<\frac{1}{p'_{\theta}}$, iff $\rho > 2$ as assumed. Furthermore 
$\frac{1}{p'_{\theta}} < \frac{1}{r_{0,\theta}}$, iff $\frac{1}{p'}<\frac{1}{q}- \theta(\frac{1}{\rho_0}-\frac{1}{2})$ (iii) and finally
$\frac{1}{p'_{\theta}} < \frac{2}{r'_{0,\theta}}$, iff $\frac{1}{p}<1 + \theta (\frac{2}{\rho}-\frac{1}{2})$ (iv). Now all the conditions
(i) - (iv) can be satisfied by choosing $\theta$ close enough to zero.

Concerning the $p=q$ - case, we observe at first that it is contained in the preceeding as long as $p=q<2$. Next we integrate (\ref{200}) with respect
to $\tau$ and $\xi$ to obtain
\[\n{uvw}{\widehat{L^{\infty}_{xt}}} \le c \n{u_0}{\widehat{L^{\infty}_x}}\n{v_0}{\widehat{L^{\infty}_x}}\n{w_0}{\widehat{L^{\infty}_x}}.\]
Finally, the claimed estimate follows for arbitrary $p=q>1$ by interpolation.
$\hfill \Box$

\vspace{0.5cm}

The $\x$-version of the above Corollary reads as follows:

\begin{kor}\label{k22} For $p$ and $q$ as in Corollary \ref{k21} the estimate
\[\n{uv\overline{w}}{\widehat{L^q_x}(\widehat{L^p_t})} \le c \n{u}{\X{q}{0}{\frac{1}{q}+}}\n{v}{\X{p}{0}{\frac{1}{p}+}}\n{w}{\X{p}{0}{\frac{1}{p}+}}.\]
is valid.
\end{kor}

\section{Local well-posedness results}

Setting $r:=p=q$ in Corollary \ref{k22}, we see that for any $b' \le 0$ the estimate
\begin{equation}\label{298}
\n{uv\overline{w}}{\X{r}{0}{b'}} \le c \n{u}{\X{r}{0}{b}}\n{v}{\X{r}{0}{b}}\n{w}{\X{r}{0}{b}}
\end{equation}
holds true. Since
\begin{equation}\label{299} 
\langle \xi \rangle \le c \sum_{i=1}^3\langle \xi_i \rangle \le c \prod_{i=1}^3\langle \xi_i \rangle,
\end{equation}
whenever $\xi = \sum_{i=1}^3 \xi_i$, we may, for $s \ge 0$, replace $\X{r}{0}{b'}$ by $\X{r}{s}{b'}$ and $\X{r}{0}{b}$ by $\X{r}{s}{b}$ in (\ref{298}).
Inserted in Theorem \ref{t0} this yields Theorem \ref{t1}. To prove Theorem \ref{t2}, some more work has to be done:

\subsection{Nonlinear estimates proving Theorem \ref{t2}} 

\begin{lemma}\label{l30} Let $r>1$, $s\ge \frac{1}{2}$, $b>\frac{1}{r}$ and $b' \le-\frac{1}{2r'}$. Then
\[\n{u_1u_2 \partial_x\overline{u}_3}{\X{r}{s}{b'}} \le c \prod_{i=1}^3 \n{u_i}{\X{r}{s}{b}}.\]
\end{lemma}

Proof:. Let $\xi_i$ denote the frequencies belonging to the $u_i$, $1 \le i \le 3$. By (\ref{299}) we may restrict ourselves to $s=\frac{1}{2}$.
Furthermore, by symmetry between the first two factors we may assume that $|\xi_1| \le |\xi_2|$. Now we consider two cases, where in the first one we suppose that 
\begin{equation}\label{300}
J^{\frac{1}{2}}(u_1u_2 \partial_x\overline{u}_3)\preceq (J^{\frac{1}{2}}u_1)(J^{\frac{1}{2}}u_2)(J^{\frac{1}{2}}\overline{u}_3).
\end{equation}
(Here and in the sequel $f \preceq g$ stands for $|\widehat{f}|\le c |\widehat{g}|$, $J$ denotes the Bessel potential operator of order $-1$.) By Corollary \ref{k22} we have
\[\n{(J^{\frac{1}{2}}u_1)(J^{\frac{1}{2}}u_2)(J^{\frac{1}{2}}\overline{u}_3)}{\widehat{L^r_{xt}}}\le c \prod_{i=1}^3 \n{u_i}{\X{r}{\frac{1}{2}}{b}},\]
as desired. Observe that (\ref{300}) holds, if $|\xi| \le 1$, or $|\xi_3| \le 1$, or, most important, if $|\xi_3| \langle \xi \rangle \le c\langle \xi_1 \rangle\langle \xi_2 \rangle$.
So, in the remaining case 2, where (\ref{300}) does not hold, we have $|\xi_1\xi_2| \ll |\xi\xi_3|$ and whence
\[\sum_{i=0}^3\langle \sigma_i \rangle \ge c |\xi_1\xi_2 - \xi\xi_3|  \ge c |\xi\xi_3| \ge c  \langle \xi \rangle \langle \xi_3 \rangle,\]
where $\sigma_0 = \tau + \xi ^2$, $\sigma_{1,2} = \tau_{1,2} + \xi_{1,2} ^2$, $\sigma_3 = \tau_3 - \xi_3 ^2$ and $\sum_{i=1}^3(\tau_i , \xi_i) =(\tau, \xi)$.
Next, we discuss the four subcases according to which one of the $\sigma$'s is the largest: 

\vspace{0,5cm}

{\bf {\underline{Subcase 2.0:}}} $\langle \sigma_0 \rangle \ge \langle \sigma_i \rangle$, $1 \le i \le 3$.\\

Suppose first that, in addition, $|\xi_3| \le c |\xi_2|$. Then $|\xi_1| \ll |\xi| \sim |\xi_2 + \xi_3|$ and $|\xi \xi_3| \le c |\xi_2(\xi_2 + \xi_3)|$, so that
\[J^{\frac{1}{2}}(u_1u_2 \partial_x\overline{u}_3)\preceq (J^{\frac{1}{2}-\frac{1}{r}}u_1)I^{\frac{1}{r}}((J^{\frac{1}{2}}u_2)(J^{\frac{1}{2}}\overline{u}_3)).\]
Cancelling another $\langle \xi_1 \rangle ^{\varepsilon}$ by $\langle \sigma_0 \rangle ^{-b'}$ we see now that
\begin{eqnarray*}
&& \n{(J^{\frac{1}{2}-\frac{1}{r}}u_1)I^{\frac{1}{r}}((J^{\frac{1}{2}}u_2)(J^{\frac{1}{2}}\overline{u}_3))}{\X{r}{0}{b'}}\\
& \le & c \n{(J^{\frac{1}{2}-\frac{1}{r}-}u_1)I^{\frac{1}{r}}((J^{\frac{1}{2}}u_2)(J^{\frac{1}{2}}\overline{u}_3))}{\widehat{L^r_{xt}}}\\
& \le & c \n{J^{\frac{1}{2}-\frac{1}{r}-}u_1}{\widehat{L^{\infty}_{xt}}}\n{I^{\frac{1}{r}}((J^{\frac{1}{2}}u_2)(J^{\frac{1}{2}}\overline{u}_3))}{\widehat{L^r_{xt}}} \le c \prod_{i=1}^3 \n{u_i}{\X{r}{\frac{1}{2}}{b}},
\end{eqnarray*}
where in the last step we have used the embedding 
\begin{equation}\label{301}
\X{r}{\frac{1}{r}+}{\frac{1}{r}+} \subset \widehat{L^{\infty}_{xt}}
\end{equation}
for the first and Corollary \ref{k20} for the second factor.

Assume now that $|\xi_2| \ll |\xi_3|$. Then $|\xi_3| \sim |\xi| \sim |\xi_2 + \xi_3|$ and thus
\begin{eqnarray*}
J^{\frac{1}{2}}(u_1u_2 \partial_x\overline{u}_3) & \preceq & I^{\frac{1}{r'}}((J^{-\frac{1}{2}}u_1)I^{\frac{1}{r}}((J^{\frac{1}{2}}u_2)(J^{\frac{1}{2}}\overline{u}_3))) \\
& \preceq & \Lambda_0^{-b'} ((J^{-\frac{1}{2}}u_1)I^{\frac{1}{r}}((J^{\frac{1}{2}}u_2)(J^{\frac{1}{2}}\overline{u}_3))) 
\end{eqnarray*}
($\Lambda_0^{-b'} = \mathcal{F}^{-1} \langle \sigma_0 \rangle^{-b'} \mathcal{F}$), the latter, since $|\xi|^{\frac{1}{r'}} \le c \langle \sigma_0 \rangle^{-b'}$ for $b' \le - \frac{1}{2r'}$.
Now
\begin{eqnarray*}
&& \n{\Lambda_0^{-b'} ((J^{-\frac{1}{2}}u_1)I^{\frac{1}{r}}((J^{\frac{1}{2}}u_2)(J^{\frac{1}{2}}\overline{u}_3)))}{\X{r}{0}{b'}}\\
& = & \n{(J^{-\frac{1}{2}}u_1)I^{\frac{1}{r}}((J^{\frac{1}{2}}u_2)(J^{\frac{1}{2}}\overline{u}_3))}{\widehat{L^r_{xt}}}\\
& \le & c \n{J^{-\frac{1}{2}}u_1}{\widehat{L^{\infty}_{xt}}}\n{I^{\frac{1}{r}}((J^{\frac{1}{2}}u_2)(J^{\frac{1}{2}}\overline{u}_3))}{\widehat{L^r_{xt}}} \le c \prod_{i=1}^3 \n{u_i}{\X{r}{\frac{1}{2}}{b}}
\end{eqnarray*}
by (\ref{301}) and Corollary \ref{k20} again. This concludes the discussion concerning subregion 2.0.

\vspace{0,5cm}

{\bf {\underline{Subcase 2.1:}}} $\langle \sigma_1 \rangle = \max_{i=0}^3 \langle \sigma_i \rangle$. \\

In this case we have
\[J^{\frac{1}{2}}(u_1u_2 \partial_x\overline{u}_3)  \preceq (J^{\frac{1}{2}-b}\Lambda_1^bu_1)(J^{\frac{1}{2}-b}u_2)(J^{\frac{1}{2}}\overline{u}_3))\]
($\Lambda_1^{b} = \mathcal{F}^{-1} \langle \sigma_1 \rangle^{b} \mathcal{F}$), where the $\widehat{L^r_{xt}}$-norm of the latter is bounded by
\[c\n{J^{\frac{1}{2}-b}\Lambda_1^bu_1}{\widehat{L^{\infty}_x}(\widehat{L^{r}_{t}})}\n{J^{\frac{1}{2}-b}u_2}{\widehat{L^{\infty}_{xt}}}\n{J^{\frac{1}{2}}u_3}{\widehat{L^{r}_{x}}(\widehat{L^{\infty}_{t})}}\le c \prod_{i=1}^3 \n{u_i}{\X{r}{\frac{1}{2}}{b}},\]
where we have used the embeddings $\X{r}{\frac{1}{r}+}{0} \subset \widehat{L^{\infty}_x}(\widehat{L^{r}_{t}})$, (\ref{301}) and $\X{r}{0}{\frac{1}{r}+} \subset \widehat{L^{r}_{x}}(\widehat{L^{\infty}_t})$.
\vspace{0,5cm}

{\bf {\underline{Subcase 2.2:}}} $\langle \sigma_2 \rangle = \max_{i=0}^3 \langle \sigma_i \rangle$. \\

- can be treated in exactly the same manner.
\vspace{0,5cm}

{\bf {\underline{Subcase 2.3:}}} $\langle \sigma_3 \rangle = \max_{i=0}^3 \langle \sigma_i \rangle$. \\

Here
\[J^{\frac{1}{2}}(u_1u_2 \partial_x\overline{u}_3)  \preceq (J^{\frac{1}{2}-b}u_1)(J^{\frac{1}{2}-b}u_2)(J^{\frac{1}{2}}\Lambda_3^b\overline{u}_3)),\]
with $\Lambda_3^{b} = \mathcal{F}^{-1} \langle \sigma_3 \rangle^{b} \mathcal{F}$, so that $\n{J^{\frac{1}{2}}\Lambda_3^b\overline{u}}{\widehat{L^r_{xt}}}= \n{u}{\X{r}{\frac{1}{2}}{b}}$.
Putting the first two factors into $\widehat{L^{\infty}_{xt}}$ and the third one into $\widehat{L^r_{xt}}$ we end up with the desired bound, after having used (\ref{301}) for the first two factors.
$\hfill \Box$

\begin{lemma}\label{l31} Let $r>1$, $s \ge \frac{1}{2}$, $b>\frac{1}{r}$. Then
\[\n{J^{s}(u_1u_2u_3 \overline{u}_4\overline{u}_5)}{\widehat{L^r_{xt}}}\le c \prod_{i=1}^5 \n{u_i}{\X{r}{s}{b}}.\]
\end{lemma}

Proof: Again, we may restrict ourselves to the case $s=\frac{1}{2}$. Let $\xi_i$ denote the frequencies belonging to the $u_i$, $1 \le i \le 5$. Without loss we may assume $|\xi_1| \ge |\xi_2| \ge |\xi_3|$
and $|\xi_4| \ge |\xi_5|$. Then, if $|\xi_1| \ge |\xi_4|$, we have
\[J^{\frac{1}{2}}(u_1u_2u_3 \overline{u}_4\overline{u}_5) \preceq (J^{\frac{1}{2}}u_1)(J^{\frac{1}{2}}u_2)(J^{\frac{1}{2}}\overline{u}_4)(J^{-\frac{1}{2}}u_3)(J^{-\frac{1}{2}}\overline{u}_5),\]
which gives the upper bound
\[\n{(J^{\frac{1}{2}}u_1)(J^{\frac{1}{2}}u_2)(J^{\frac{1}{2}}\overline{u}_4)}{\widehat{L^r_{xt}}}\n{J^{-\frac{1}{2}}u_3}{\widehat{L^{\infty}_{xt}}}\n{J^{-\frac{1}{2}}u_3}{\widehat{L^{\infty}_{xt}}} \le c  \prod_{i=1}^5 \n{u_i}{\X{r}{\frac{1}{2}}{b}}\]
by Corollary \ref{k22} (for the first factor) and by the embedding (\ref{301}) (for the last two factors).
If $|\xi_1| \le |\xi_4|$, we consider two subcases; either $|\xi_5| \le |\xi_1|$, where we get the same bound as above, or $|\xi_5| \ge |\xi_1|$, where we have
\[J^{\frac{1}{2}}(u_1u_2u_3 \overline{u}_4\overline{u}_5) \preceq (J^{\frac{1}{2}}u_1)(J^{\frac{1}{2}}\overline{u}_4)(J^{\frac{1}{2}}\overline{u}_5)(J^{-\frac{1}{2}}u_2)(J^{-\frac{1}{2}}u_3).\]
Now Corollary \ref{k22} (observe that \n{fg\overline{h}}{\widehat{L^r_{xt}}}=\n{\overline{fg}h}{\widehat{L^r_{xt}}}) and (\ref{301})
again lead to the desired bound.
$\hfill \Box$

\subsection{Continuity of the gauge transform in $\h{r}{s}$ - spaces and proof of Theorem \ref{t3}} 

\begin{lemma}\label{l32}
Let $s \ge \frac{1}{2}$, $2\ge r>1$. Then $G$ as a map from $\h{r}{s}$ to $\h{r}{s}$ (and from $C(I,\h{r}{s})$ to $C(I,\h{r}{s})$, respectively)
is Lipschitz contiuous on bounded subsets.
\end{lemma}

Proof: If $m: \R \rightarrow \C$ is of bounded variation over $\R$ and $M$ is the Fourier multiplier associated with $m$, then $M$ is a bounded
operator from $L^p(\R)$ to $L^p(\R)$, $1<p<\infty$, and for the operator norm there is the bound
\[\|M\|_{L^p \rightarrow L^p} \le c \,\,\,\,\,\,(\lim_{x \rightarrow -\infty} |m(x)| + \int_{-\infty}^{\infty}|dm(x)|\,\,\,)\]
(see \cite[Corollary 3.8]{D01} and its proof). From the definition of the spaces $\widehat{L^r_x}$ it follows that a pointwise multiplier on
$\widehat{L^r_x}$ acts like a Fourier multiplier on $L^{r'}_{\xi}$. Hence, for $u \in \widehat{L^r_x}$,
\[\n{mu}{\widehat{L^r_x}}\le c \,\,\,\,\,\,(\lim_{x \rightarrow -\infty} |m(x)| + \int_{-\infty}^{\infty}|dm(x)|\,\,\,)\n{u}{\widehat{L^r_x}}.\]
Fixing $v \in L^2_x$ and writing $G_v(x):= e^{-i \int_{-\infty}^x |v(y)|^2 dy}$ we obtain
\begin{equation}\label{302}
\n{G_vu}{\widehat{L^r_x}} \le c (1+\q{v}{L^2_x})\n{u}{\widehat{L^r_x}}.
\end{equation}
If, in addition, $w \in L^2_x$, we use the above and the mean value theorem to see that
\begin{equation}\label{303}
\n{(G_v-G_w)u}{\widehat{L^r_x}} \le c \n{v-w}{L^2_x}(1+\n{v}{L^2_x}+\n{w}{L^2_x})^3\n{u}{\widehat{L^r_x}}.
\end{equation}
Next we consider $u \in \h{r}{1}$ and $v,w \in \h{r}{\frac{1}{2}}$. Then
\begin{eqnarray*}
\n{G_vu}{\h{r}{1}} & \le & \n{G_vu}{\widehat{L^r_x}} + \n{\partial_x(G_vu)}{\widehat{L^r_x}} \\
& \le & \n{G_vu}{\widehat{L^r_x}}+\n{G_v|v|^2u}{\widehat{L^r_x}}+\n{G_vu_x}{\widehat{L^r_x}}.
\end{eqnarray*}
Using (\ref{302}) we get
\[\n{G_vu}{\h{r}{1}}  \le c(1+\q{v}{L^2_x})(\n{u}{\h{r}{1}} + \n{|v|^2u}{\widehat{L^r_x}}).\]
Now the Hausdorff-Young inequality and the embedding $\h{r}{s} \subset \widehat{L^{\rho}_x}$, $s>\frac{1}{r}-\frac{1}{\rho}\ge 0$ give
\[\n{|v|^2u}{\widehat{L^r_x}} \le c \q{v}{\widehat{L^{2r}_x}}\n{u}{\widehat{L^{\infty}_x}} \le c \q{v}{\h{r}{\frac{1}{2}}}\n{u}{\h{r}{1}},\]
hence
\[\n{G_vu}{\h{r}{1}} \le c (1+\q{v}{L^2_x})(1+\q{v}{\h{r}{\frac{1}{2}}})\n{u}{\h{r}{1}}.\]
Interpolation with (\ref{302}) yields for $0 \le s \le 1$
\[\n{G_vu}{\h{r}{s}} \le c (1+\q{v}{L^2_x})(1+\q{v}{\h{r}{\frac{1}{2}}})\n{u}{\h{r}{s}}.\]
Similar estimates using in addition (\ref{303}) show that
\[\n{(G_v-G_w)u}{\h{r}{s}} \le c \n{v-w}{\h{r}{\frac{1}{2}}}(1+\n{v}{\h{r}{\frac{1}{2}}}+\n{w}{\h{r}{\frac{1}{2}}})^5\n{u}{\h{r}{s}},\]
where still $0\le s\le 1$. Especially, if $\frac{1}{2} \le s \le 1$ and $u,v,w \in \h{r}{s}$ we have
\[\n{G_vu}{\h{r}{s}} \le c (1+\n{v}{\h{r}{s}})^4\n{u}{\h{r}{s}}\]
and
\[\n{(G_v-G_w)u}{\h{r}{s}} \le c\n{v-w}{\h{r}{s}}(1+\n{v}{\h{r}{s}}+\n{w}{\h{r}{s}})^5\n{u}{\h{r}{s}}.\]
Concerning higher regularity we use induction and similar arguments to show that for $s>1$ there exist exponents $\alpha=\alpha(s)$ such that
\begin{equation}\label{305}
\n{G_vu}{\h{r}{s}} \le c (1+\n{v}{\h{r}{s}})^{\alpha}\n{u}{\h{r}{s}}
\end{equation}
and
\begin{equation}\label{306}
\n{(G_v-G_w)u}{\h{r}{s}} \le c\n{v-w}{\h{r}{s}}(1+\n{v}{\h{r}{s}}+\n{w}{\h{r}{s}})^{\alpha}\n{u}{\h{r}{s}}.
\end{equation}
Finally, for $u,v \in \h{r}{s}$, $s \ge \frac{1}{2}$ we obtain by (\ref{305}) and (\ref{306})
\begin{eqnarray}\label{307}
\n{Gu-Gv}{\h{r}{s}} & \le & \n{(G_u-G_v)u}{\h{r}{s}} +\n{G_v(u-v)}{\h{r}{s}} \nonumber \\
& \le & c (1+\n{u}{\h{r}{s}}+\n{v}{\h{r}{s}})^{\alpha + 1}\n{u-v}{\h{r}{s}}.
\end{eqnarray}
The proof is completed by the remark that for time dependent functions $u,v \in C(I,\h{r}{s})$ we obviously may replace $\|\,\,\,\|_{\h{r}{s}}$
by $\sup_{t\in I}\|\,\,\,\|_{\h{r}{s}}$ in (\ref{307}).
$\hfill \Box$

{\bf{Remark:}} Obviously, Lemma \ref{l32} is equally valid with $G$ replaced by $G^{-1}$, where the inverse transform $G^{-1}$ is given by
\[G^{-1}v(x):=e^{i\int_{- \infty}^x |v(y)|^2 dy}v(x).\]

Proof of Theorem \ref{t3}: If $u_0 \in \h{r}{s}$, then so is, by Lemma \ref{l32}, $v_0:=Gu_0$. Theorem \ref{t2} gives a unique solution 
\[v \in \x(\delta) \subset C([-\delta,\delta], \h{r}{s})\]
of (\ref{110}). Moreover, for $\delta_0 \in (0,\delta)$ the mapping $S:u_0 \mapsto u$, $\h{r}{s} \rightarrow \x(\delta_0)$ is locally Lipschitz
continuous. Now $u:=G^{-1}v$ solves the IVP (\ref{109}). For smooth solutions this is clear by formal computations (cf. e. g. \cite[p. 1498]{HO94}),
while the general case can be reduced to this by approximation as follows: Let $u^{(n)}_0 \in \mathcal{S}(\R)$ be a sequence of data with $u^{(n)}_0 \rightarrow u_0$
in $\h{r}{s}$ and $u^{(n)} = G^{-1}SG u^{(n)}_0$. Then $u^{(n)}\rightarrow u$ in $C([-\delta_0, \delta_0],\h{r}{s})$, $Gu^{(n)}\rightarrow Gu$ in $\x(\delta_0)$ and,
for $|t|\le \delta_0$,
\begin{equation}\label{308}
u^{(n)}(t) = e^{it\partial_x^2}u^{(n)}_0 + \int_0^t e^{i(t-s)\partial_x^2}(|u^{(n)}|^2u^{(n)})_x(s)ds.
\end{equation}
Clearly $u^{(n)}(t) \rightarrow u(t)$ and $e^{it\partial_x^2}u^{(n)}_0\rightarrow e^{it\partial_x^2}u_0$ in $\h{r}{s}$. Next we use the embeddings
$L^1_x \subset H_x^{-\frac{1}{2}-}$, $\x \subset \XX{0}{\frac{3}{8}+} \subset L^4_x $ and $\h{r}{s} \subset L^2_x$ to estimate
\[\sup_{|t|\le \delta_0}\|\int_0^t e^{i(t-s)\partial_x^2}(|u^{(n)}|^2u^{(n)}-|u|^2u)_x(s)ds\|_{H_x^{-\frac{3}{2}-}}\]
coarsely by
\[c(\q{Gu^{(n)}}{\x(\delta_0)} + \q{Gu}{\x(\delta_0)})\n{Gu^{(n)}-Gu}{L_t^{\infty}([-\delta_0, \delta_0], \h{r}{s})},\]
which tends to zero for $n \rightarrow \infty$. Hence (\ref{308}) holds with $u^{(n)}$ and $u^{(n)}_0$ replaced by $u$ and $u_0$, thus existence
is shown. Uniqueness of $u$ follows from that of $v$. Persistence property and the statement on continuous dependence are now immediate consequences
of Lemma \ref{l32}.
$\hfill \Box$

\section{Remarks on ill-posedness}

This section is devoted to review the arguments from \cite{KPV01} and \cite{BL01}, respectively, showing local ill-posedness
for cubic (focusing\footnote{In contrast to defocusing, i. e.: with the opposite sign before the nonlinearity in (\ref{101}); ill-posedness results 
concerning this case are obtained in \cite{CCT03}.}) NLS below $L^2_x$ and, respectively, for DNLS below $H_x^{\frac{1}{2}}$. By ill-posedness it is meant here
that the mapping data upon solution, even when restricted to closed balls of the data space, cannot be uniformly continuous into
any solution space being continuously embedded into the continuous functions on a time interval $[0,T]$ with values in the data space (cf. \cite{KPV01}, p. 617 f.).
It turns out that these arguments work well - with minor changes - when data in the spaces $\h{r}{s}$ are considered.

\subsection{Ill-posedness of cubic (focusing) NLS in $\h{r}{s}$ for $s<0$} \hfill \\

The following counterexample was given in \cite{KPV01} in order to show that the Cauchy-problem (\ref{101})
is locally ill posed for $u_0 \in H^s(\R)$ if $s < 0$:

\vspace{0.3cm}

Let $f(x)= \frac{\sqrt{2}}{\cosh{(x)}}$. Then $f$ solves the ODE $f'' - f +  f^3 = 0$.
Setting $f_{\omega}(x) = \omega f(\omega x)$ and
\[u_{N \omega}(x,t)=\exp{(-it(N^2 - \omega ^2) + iNx)}f_{\omega}(x- 2N t)\]
one gets a two parameter family of solutions of (\ref{101}) with data
\[u_0(x)=u_{N \omega}(x,0)=\exp{(iNx)}f_{\omega}(x).\]
Now sequences $N_{1,2} \sim N \rightarrow \infty$ and $\omega = N^{-2s}$ are chosen so that for $-\frac{1}{2}<s<0$
\[\|u_{N_1 \omega}(0) - u_{N_2 \omega}(0)\|_{H^s}\le c \omega ^{-\frac{1}{2}}N^s |N_1 -N_2| = c N^{2s}|N_1 -N_2|\]
and
\[\|u_{N_1 \omega}(T) - u_{N_2 \omega}(T)\|_{H^s}\ge c ,\]
provided $|N_1 -N_2|T \gg \omega ^{-1}$ respectively $|N_1 -N_2| \gg \frac{N^{2s}}{T}$ ensuring that the supports of the $u_{N_j \omega}(T), \,\,j=1,2$ are essentially disjoint. 
Now if $N_1 -N_2 = \frac{C}{T}N^{2s}$ with a large constant $C$, the latter condition is fulfilled and
\[\|u_{N_1 \omega}(0) - u_{N_2 \omega}(0)\|_{H^s}\le c \frac{C}{T}N^{4s} \longrightarrow 0 \,\,\,\,(s < 0).\]
Thus the mapping data upon solution from $H^s(\R)$ to any solution space $X_T$ continuously embedded in $C([0,T],H^s(\R))$ 
cannot be uniformly continuous.

\vspace{0.3cm}

So far, this is nothing but a short summary of the argument given by Kenig, Ponce and Vega, see Thm. 1.1 and \S 2 in \cite{KPV01} for the 
details. The same example shows local illposedness of (\ref{101}) for data $u_0 \in \h{r}{s}(\R)$, if $r>1$, $-\frac{1}{r'} < s < 0$. 
In fact, if we follow step by step the computations in \cite{KPV01} and choose $\omega = N^{-sr'}$, we see that - since the $u_{N_j \omega}$ are frequency concentrated around $N$ by the assumption $-\frac{1}{r'} < s$ -
\[\|u_{N_1 \omega}(0) - u_{N_2 \omega}(0)\|_{\h{r}{s}}\le c \omega ^{-\frac{1}{r}}N^s |N_1 -N_2| = c N^{sr'}|N_1 -N_2|.\]
On the other hand we have
\begin{eqnarray*}
\|u_{N_1 \omega}(T) - u_{N_2 \omega}(T)\|_{\h{r}{s}} & \ge & c N^s \|u_{N_1 \omega}(T) - u_{N_2 \omega}(T)\|_{\widehat{L^r_x}} \\ 
& \ge & c N^s \sup_{\|\phi\|_{L^r} \le 1}\langle u_{N_1 \omega}(T) - u_{N_2 \omega}(T), \check{\phi}\rangle _{L^2_x}.
\end{eqnarray*}
Now the $u_{N_j \omega}(T), \,\,j=1,2$, are concentrated on intervals $I_j$ of size $\omega ^{-1}$ around $2N_jT$, which are disjoint for (cf. (2.17) in \cite{KPV01})
\begin{equation}\label{401}
|N_1-N_2| T \gg \omega ^{-1} =N^{sr'}.
\end{equation}
Choosing
\[\check{\phi}= c \omega ^{\frac{1}{r'}} \chi_{I_1} \frac{\overline{u_{N_1 \omega}(T)}}{|u_{N_1 \omega}(T)|},\]
where the factor $\omega ^{\frac{1}{r'}}$ ensures that $\|\phi\|_{L^r} \le 1$,
we obtain the lower bound
\[c N^s \omega ^{\frac{1}{r'}} \int |u_{N_1 \omega}(x,T)|dx = c \omega \int f(\omega x) dx = c .\]
With $N_1 - N_2 = \frac{C}{T} N ^ {sr'}$ (so that (\ref{401}) is fulfilled) we have
\[\|u_{N_1 \omega}(0) - u_{N_2 \omega}(0)\|_{\h{r}{s}} \le c N^{2sr'},\]
which tends to zero, if $s < 0$ and $r>1$.

\subsection{Ill-posedness of DNLS in $\h{r}{s}$ for $s<\frac{1}{2}$} \hfill \\

Concerning the Cauchy problem (\ref{109}) we can rely on the work of Biagioni and Linares, see Theorem 2.1 in \cite{BL01} and its proof.
There the following two parameter family of solutions of DNLS was used to build up a counterexample showing ill-posedness of (\ref{109}) in $H^s$, if $s<\frac{1}{2}$:
Let
\begin{eqnarray*}
N \ge 0, \hspace{0.5cm} \omega > \frac{N^2}{4},\hspace{0.5cm} \gamma^2 = 4\omega - N^2, \hspace{0.5cm} \alpha = \frac{N}{2\sqrt{\omega}},\hspace{0.5cm} \mbox{and}\\
\phi(x)= 3 \arctan{\left(\frac{\exp{(x)+\alpha}}{\sqrt{1- \alpha^2}}\right)},\hspace{0.3cm} f(x)=(\cosh{(x)}+\alpha)^{-\frac{1}{2}},\\
F(x)=e^{i\phi(x)}f(x),\hspace{2cm} F_{\gamma}(x)=\gamma F(\gamma x).
\end{eqnarray*}
Then $u_{N \omega}$, defined by
\[u_{N \omega}(x,t)=\exp{\left(i\left(\frac{Nx}{2}+(\omega- \frac{N^2}{2})t\right)\right)}\omega^{-\frac{1}{4}}F_{\gamma}(x-Nt)\]
solves (\ref{109}) with data
\[u_0(x)=u_{N \omega}(x,0)=\exp{\left(i\frac{Nx}{2}\right)}\omega^{-\frac{1}{4}}F_{\gamma}(x).\]
Cf. (2.2) in \cite{BL01}, here we use a slightly different notation. This family of solutions of DNLS has been derived by Hohenberg and van Saarloos
in a more general context, see section 3.2.3.1 in \cite{vSH92}. It shall be used in the sequel to generalize the result of Biagioni
and Linares to data in $\h{r}{s}$, when $r>1$ and $\frac{1}{2} > s > \frac{1}{2}-\frac{1}{r'}$. For that purpose we choose, for given $N\ge 0$,
\hfill \\
\begin{itemize}
\item $4\omega = N^2 + N^{r'(1-2s)}$, $\gamma = N^{r'(\frac{1}{2}-s)}$, $\alpha = \frac{N}{2\sqrt{\omega}}$,
\hfill \\
\item $N'=N+C$ with a large positive constant $C$,
\hfill \\
\item $4\omega' = {N'} ^2 + N^{r'(1-2s)}\frac{{N'}^2}{N^2}$, $\gamma'=\gamma \frac{N'}{N}$, $\alpha' = \frac{N'}{2\sqrt{\omega'}}$.
\end{itemize}
\hfill \\
Then $\gamma^2=4\omega - N^2$, ${\gamma'}^2=4\omega' - {N'}^2$ and $\alpha^2={\alpha'}^2<1$. The condition $s > \frac{1}{2}-\frac{1}{r'}$ 
assures that $\omega \sim N^2$ for large $N$ and $\frac{N}{\gamma}\ge1$, similarly for $N'$, $\omega'$ and $\gamma'$, which shall be used
in the subsequent computations. Now it is sufficient to show that
\begin{equation}\label{402}
\n{u_{N \omega}(0)-u_{N' \omega'}(0)}{\h{r}{s}} \longrightarrow 0\hspace{0,5cm}(N \longrightarrow \infty)
\end{equation}
while
\begin{equation}\label{403}
\n{u_{N \omega}(T)-u_{N' \omega'}(T)}{\h{r}{s}} \ge c.
\end{equation}

To see (\ref{402}) we observe first that
\begin{eqnarray*}
&&|\widehat{u_{N \omega}(0)}(\xi)-\widehat{u_{N' \omega'}(0)}(\xi)|=|\omega^{-\frac{1}{4}}\widehat{F}(\frac{\xi - \frac{N}{2}}{\gamma})-{\omega'}^{-\frac{1}{4}}\widehat{F}(\frac{\xi - \frac{N'}{2}}{\gamma'})|\\
& \le & \omega^{-\frac{1}{4}}|\widehat{F}(\frac{\xi - \frac{N}{2}}{\gamma})-\widehat{F}(\frac{\xi - \frac{N'}{2}}{\gamma'})|+|\omega^{-\frac{1}{4}}-{\omega'}^{-\frac{1}{4}}||\widehat{F}(\frac{\xi - \frac{N'}{2}}{\gamma'})|=:I+II
\end{eqnarray*}
with
\[\|I\|^{r'}_{\h{r}{s}}=\omega^{-\frac{r'}{4}}\int\langle \xi\rangle^{r's}|\widehat{F}(\frac{\xi - \frac{N}{2}}{\gamma})-\widehat{F}(\frac{\xi - \frac{N'}{2}}{\gamma'})|^{r'}d\xi.\]
Writing $G(x)=\widehat{F}(x-\frac{N}{2\gamma})$ ($=\widehat{F}(x-\frac{N'}{2\gamma'})$ by our choice of parameters) and substituting $\eta = \frac{\xi}{\gamma}$ we have
\[\|I\|^{r'}_{\h{r}{s}} \le c N^{-\frac{r'}{2}}\gamma^{r's+1}\int \langle \eta\rangle^{r's}|G(\eta)-G(\eta \frac{\gamma}{\gamma'})|^{r'}d \eta\]
with
\[|G(\eta)-G(\eta \frac{\gamma}{\gamma'})|^{r'} =|\int^{\eta}_{\eta \frac{\gamma}{\gamma'}}G'(\xi)d\xi|^{r'}\le (|\eta|(1-\frac{\gamma}{\gamma'}))^{\frac{r'}{r}}\int^{\eta}_{\eta \frac{\gamma}{\gamma'}}|G'(\xi)|^{r'}d\xi\]
by the mean value theorem and H\"older's inequality. This gives
\begin{eqnarray*}
&&\int \langle \eta\rangle^{r's}|G(\eta)-G(\eta \frac{\gamma}{\gamma'})|^{r'}d \eta \\
& \le & \left(\frac{C}{N}\right)^{\frac{r'}{r}}\int \langle \eta\rangle^{r'(s+\frac{1}{r})}\int^{\eta}_{\eta \frac{\gamma}{\gamma'}}|G'(\xi)|^{r'}d\xi d \eta \\
& = & \left(\frac{C}{N}\right)^{\frac{r'}{r}}\int|G'(\xi)|^{r'}\int_{\xi}^{\xi\frac{\gamma'}{\gamma}}\langle \eta\rangle^{r'(s+\frac{1}{r})}d \eta d\xi\\
& \le & c \left(\frac{C}{N}\right)^{\frac{r'}{r}} \left(\frac{\gamma'}{\gamma}-1\right) \int \langle \xi \rangle^{r'(s+1)}|G'(\xi)|^{r'}d\xi \\
&=& c \left(\frac{C}{N}\right)^{r'}\int \langle \xi \rangle^{r'(s+1)}|G'(\xi)|^{r'}d\xi.
\end{eqnarray*}
Now
\begin{eqnarray*}
\int \langle \xi \rangle^{r'(s+1)}|G'(\xi)|^{r'}d\xi &=& \int \langle \xi \rangle^{r'(s+1)}|\widehat{F}'(\xi-\frac{N}{2\gamma})|^{r'}d\xi \\
& \le & c \left(\frac{N}{\gamma}\right)^{r'(s+1)}\int \langle \xi \rangle^{r'(s+1)}|\widehat{F}'(\xi)|^{r'}d\xi,
\end{eqnarray*}
where the last integral, although dependent on the parameter $\alpha \in (0,1)$, is bounded by a constant. We arrive at
\[\n{I}{\h{r}{s}} \le c N^{r'(s-\frac{1}{2})} \longrightarrow 0\,\,\,\,\,(N\longrightarrow \infty).\]
In order to estimate $\n{II}{\h{r}{s}}$ we first notice that $|\omega^{-\frac{1}{4}}-{\omega'}^{-\frac{1}{4}}|\sim N^{-\frac{3}{2}}$, while
\begin{eqnarray*}
\int \langle \xi \rangle^{r's}|\widehat{F}(\frac{\xi - \frac{N'}{2}}{\gamma'})|^{r'}d\xi &\le & c {\gamma'}^{r's+1}\int \langle \xi \rangle^{r's}|\widehat{F}(\xi - \frac{N'}{2\gamma'})|^{r'}d\xi \\
& \le & c \gamma' {N'}^{r's}\int \langle \xi \rangle^{r's}|\widehat{F}(\xi )|^{r'}d\xi \le c {N'}^{\frac{r'}{2}}.
\end{eqnarray*}
This gives
\[\n{II}{\h{r}{s}} \le c N^{-1} \longrightarrow 0\,\,\,\,\,(N\longrightarrow \infty).\]
Thus (\ref{402}) is shown. To obtain (\ref{403}) we define
\[\psi(x)=N^s\gamma^{\frac{1}{r'}}\exp{\left(-\frac{\gamma ^2}{2}(x-NT)^2+ i((\omega-\frac{N^2}{2})T+\frac{Nx}{2})\right)}.\]
Then
\[|\widehat{\psi}(\xi)|=cN^s\gamma^{-\frac{1}{r}}\exp{\left(-\frac{1}{2 \gamma^2}(\xi-\frac{N}{2})^2\right)}\]
and
\[\|\psi\|^r_{\h{r'}{-s}}=cN^{sr}\gamma^{-1}\int\langle \xi \rangle^{-rs}\exp{\left(-\frac{r}{2 \gamma^2}(\xi-\frac{N}{2})^2\right)}d\xi \le c.\]
Hence
\[\n{u_{N \omega}(T)-u_{N' \omega'}(T)}{\h{r}{s}} \ge c|\langle u_{N \omega}(T)-u_{N' \omega'}(T), \psi \rangle_{L^2_x}|.\]
Now $u_{N \omega}(T)$ respectively $u_{N' \omega'}(T)$ is concentrated in $[NT-\gamma^{-1},NT+\gamma^{-1}]$ respectively in $[N'T-{\gamma'}^{-1},N'T+{\gamma'}^{-1}]$,
which are disjoint for $(N'-N)T \gg \max{(\gamma^{-1},{\gamma'}^{-1})}$ (cf. (2.18) in \cite{BL01}). The latter is guaranteed by our choice of parameters and we get the lower bound
\begin{eqnarray*}
&& c |\langle u_{N \omega}(T), \psi \rangle_{L^2_x}|\\
& \ge & c N^{s-\frac{1}{2}}\gamma^{\frac{1}{r'}+1}\left|\int f(\gamma(x-NT))\exp{\left(i \phi(\gamma(x-NT)) -\frac{\gamma ^2}{2}(x-NT)^2\right)}dx\right|\\
&=& c N^{s-\frac{1}{2}}\gamma^{\frac{1}{r'}}\left|\int f(y)\exp{(i\phi(y)-\frac{y^2}{2})}dy\right| \ge c \left|\int_{-1}^1 f(y)\exp{(i\phi(y))}dy \right|.
\end{eqnarray*}
For large $N$ we have that $1-\alpha^2 \ll 1$, whence $\phi(y) \approx \frac{3\pi}{2}$, so that $|\langle u_{N \omega}(T), \psi \rangle_{L^2_x}|\ge \nolinebreak c$.
Now (\ref{403}) is established, too.

\section{Global well-posedness for cubic NLS in $\widehat{L^r_x}$, $2>r>\frac{5}{3}$, - \\ Proof of Theorem \ref{t4}}

\subsection{Splitting of the data} \hfill \\

We decompose $u_0 = u_{\le} + u_>$ with $\widehat{u_{\le}} (\xi) = \widehat{u_0} (\xi) \chi_{\{\widehat{u_0} \le \frac{1}{N}\}}(\xi)$.
Then $\widehat{u_{\le}} \in \widehat{L^{\rho}_x}$ for all $\rho \in [1,r]$ and, by convexity,
\begin{equation}\label{500}
\n{u_{\le}}{\widehat{L^{\rho}_x}} \le \|u_0\|^{\frac{r'}{\rho '}}_{\widehat{L^r}_x} N^{\frac{r'}{\rho '}-1}.
\end{equation}
On the other hand we have $u_> \in L^2$ with
\[\n{u_>}{L^2_x} \le \n{\widehat{u_0}}{L^{r'}_{\xi}} \n{\chi_{\{\widehat{u_0}\}>\frac{1}{N}}}{L^q_{\xi}} \hspace{1cm}(\frac{1}{2} = \frac{1}{r'}+\frac{1}{q}),\]
where
\[\|\chi_{\{\widehat{u_0}>\frac{1}{N}\}}\|^q_{L^q_{\xi}} = \lambda (\{\widehat{u_0}>\frac{1}{N}\}) \le N^{r'}\|\widehat{u_0}\|^{r'}_{L^{r'}_{\xi}},\]
by Tschebychev's inequality (here $\lambda$ denotes the Lebesgue-measure), which gives
\begin{equation}\label{501}
\n{u_>}{L^2_x} \le \|u_0\|^{\frac{r'}{2}}_{\widehat{L^r_x}}N^{\frac{r'}{2}-1}.
\end{equation}
We want to obtain a solution of the Cauchy-problem (\ref{101}) in the form $u=v+w$, where $v$ solves
\begin{equation}\label{502}
iv_t +v_{xx} + |v|^2 v=0 \hspace{3cm}v(0)=u_> \in L^2_x.
\end{equation}
By Strichartz' estimate ($\XX{0}{\frac{1}{2}+} \subset L^6_{xt}$) one gets a local solution $v \in \XX{0}{\frac{1}{2}+}(\delta)$ of (\ref{502}), which can be extended globally by the $L^2$-norm-conservation.
The stepwidth $\delta$ in this extension process is at most
\begin{equation}\label{503}
\delta = c\|u_>\|^{-4-}_{L^2_x} \ge c N^{(4-2r')-},
\end{equation}
since
\[\|\Lambda v\|_{\XX{0}{\frac{1}{2}+}(\delta)} \le c \n{u_>}{L^2_x} +c \delta ^{\frac{1}{2}-}\| v\|^3_{\XX{0}{\frac{1}{2}+}(\delta)},\]
$\Lambda $ being the map corresponding to the integral equation equivalent with (\ref{502}). For this solution $v$ we have - from its construction
via the contraction mapping principle - the bound
\begin{equation}\label{504}
\| v\|_{\XX{0}{\frac{1}{2}+}(\delta)} \le c \n{u_>}{L^2_x} \le c N^{\frac{r'}{2}-1}.
\end{equation}

\subsection{Local solutions for the difference equation}\hfill \\

Next we want to show existence and uniqueness of a (local, at first) solution of the Cauchy-problem for the difference equation
\begin{equation}\label{505}
iw_t + w_{xx} + N(v,w)=0, \hspace{3cm} w(0)= u_{\le}\in \widehat{L^{\rho}_x},
\end{equation}
(where $v$ is as above and
\[N(v,w)= 2|v|^2w + v^2 \overline{w} + 2v|w|^2 + w^2\overline{v} + |w|^2w)\]
with \emph{the same} (up to a constant depending on $\rho$ but not on $N$) \emph{lifespan $\delta$ as in (\ref{503})} and with - for $1 < \rho \le r$ -
\begin{equation}\label{506}
\|w\|_{\X{\rho}{0}{\frac{1}{\rho}+}(\delta)} \le c \n{u_{\le}}{\widehat{L^{\rho}_x}} \le c N^{\frac{r'}{\rho '}-1}.
\end{equation}
For this purpose, further estimates are needed:
\begin{lemma}\label{l50} For the expression $\n{fg\overline{h}}{\widehat{L^{\rho}_{xt}}}$ we have the following upper bounds:
\begin{itemize}
\item[i)] $c \n{f}{\XX{0}{\frac{1}{2}+}} \n{g}{\XX{0}{\frac{1}{2}+}} \n{h}{\X{\rho}{0}{\frac{1}{2}+}}$, if $2 \ge \rho > 1$,
\item[ii)] $c \n{f}{\X{\rho}{0}{\frac{1}{\rho}+}} \n{g}{\X{\rho}{0}{\frac{1}{\rho}+}} \n{h}{\XX{0}{b}}$, if $2 \ge \rho > \frac{4}{3}$, $b > \frac{1}{2 \rho '} + \frac{1}{4}$,
\item[iii)] $c \n{f}{\XX{0}{b}} \n{g}{\X{\rho_0}{0}{\frac{1}{\rho_0}+}} \n{h}{\X{\rho_0}{0}{\frac{1}{\rho_0}+}}$, if $\rho_0 > \frac{4}{3} \ge \rho >1$, $b> \frac{3}{2 \rho '}+\frac{1}{\rho _0}-\frac{3}{4} \ge 0$.
\end{itemize}
In any of these estimates $f$ and $h$ may be interchanged.
\end{lemma}

Proof: Part i) follows by interpolation between Corollary \ref{k22} and 
\[\n{fg\overline{h}}{\widehat{L^{\rho}_{xt}}} \le c \n{f}{\widehat{L^{\infty}_{xt}}} \n{g}{\widehat{L^{\infty}_{xt}}} \n{h}{\widehat{L^{\rho}_{xt}}}\hspace{1cm}\mbox{(Young)}.\]
To prove ii) we first use the Hausdorff-Young- and H\"older-inequalities to get
\[\n{fg\overline{h}}{\widehat{L^{\rho}_{xt}}} \le c\n{fg\overline{h}}{L^{\rho}_{xt}} \le c \n{f}{L^{3 \rho}_{xt}} \n{g}{L^{3 \rho}_{xt}} \n{h}{L^{3 \rho}_{xt}}.\]
By the $\x$-version of the Fefferman-Stein-estimate (\ref{102}) we have $\n{f}{L^{3 \rho}_{xt}} \le c \n{f}{\X{\rho}{0}{\frac{1}{\rho}+}}$ as well as
$\n{g}{L^{3 \rho}_{xt}} \le c \n{g}{\X{\rho}{0}{\frac{1}{\rho}+}}$, while $\n{h}{L^{3 \rho}_{xt}} \le c \n{h}{\XX{0}{b}}$ with $b$ as demanded follows by interpolation between
$\XX{0}{\frac{1}{2}+} \subset L^6_{xt}$ (Strichartz) and $\XX{0}{0}=L^2_{xt}$. The proof of iii) follows the same lines and will therefore be omitted.
$\hfill \Box$

In order to extract a positive power of $\delta$ from the nonlinear estimates we shall use:
\begin{lemma}\label{l51}
Let $0 < \delta \le 1$, $1 < r < \infty$, $\frac{1}{r} > b > b' \ge 0$ or $0 \ge b > b' > -\frac{1}{r'}$. Then
\[\n{f}{\X{r}{0}{b'}(\delta)} \le c \delta ^{b-b'-} \n{f}{\X{r}{0}{b}(\delta)}\]
\end{lemma}

Proof:  Let $\psi$ be a smooth cut-off-function with $\psi|_{[-1,1]}=1 $ and $\psi_{\delta}(t)=\psi(\frac{t}{\delta})$. 
Then the claimed estimate will follow from
\begin{equation}\label{507}
\n{\psi_{\delta}f}{\X{r}{0}{b'}}\le c \delta ^{b-b'-} \n{f}{\X{r}{0}{b}}.
\end{equation}
Here, by duality, it is sufficient to consider the case where $\frac{1}{r} > b > b' \ge 0$. Now (\ref{507}) can be deduced from
\begin{equation}\label{508}
\n{J_t^{b'}(\psi_{\delta}g)}{\widehat{L^r_t}}\le c \delta ^{b-b'-} \n{J_t^{b}g}{\widehat{L^r_t}}
\end{equation}
by taking $g=e^{-it\partial ^2}f$ and integrating with respect to the $\xi$-variable. Now, for $\frac{1}{r}=\frac{1}{r_1}+\frac{1}{r_2}=\frac{1}{\rho_1}+\frac{1}{\rho_2}$,
\[\n{J_t^{b'}(\psi_{\delta}g)}{\widehat{L^r_t}}\le c (\n{I_t^{b'}\psi_{\delta}}{\widehat{L^{r_1}_t}}\n{g}{\widehat{L^{r_2}_t}}+
\n{\psi_{\delta}}{\widehat{L^{\rho_1}_t}}\n{J^{b'}g}{\widehat{L^{\rho_2}_t}} ),\]
where $\n{I_t^{b'}\psi_{\delta}}{\widehat{L^{r_1}_t}}=c \delta ^{\frac{1}{r_1}-b'}\n{\psi}{\widehat{L^{r_1}_t}}$ and 
$\n{\psi_{\delta}}{\widehat{L^{\rho_1}_t}}=c \delta ^{\frac{1}{\rho_1}}\n{\psi}{\widehat{L^{\rho_1}_t}}$. Choosing $\frac{1}{r_1}=b-$, $\frac{1}{\rho_1}=b-b'-$,
we get the upper bound
\[...\le c \delta ^{b-b'-}(\n{g}{\widehat{L^{r_2}_t}}+\n{J_t^{b'}g}{\widehat{L^{\rho_2}_t}} )\]
with $\frac{1}{r_2}=(\frac{1}{r}-b)+$ and $\frac{1}{\rho_2}=(\frac{1}{r}-b+b')+$, so that $b-\frac{1}{r}> \max{(-\frac{1}{r_2}, b'-\frac{1}{\rho_2})}$. Finally,
(\ref{508}) follows by a simple H\"older-application.
$\hfill \Box$

\vspace{0,5cm}

Now let us turn to the Cauchy-problem (\ref{505}), respectively to the integral equation corresponding to it, i. e.:
\[w(t)=\Lambda_v w (t):= e^{it\partial^2}u_{\le} + i \int_0^t e^{i(t-s)\partial_x^2} N(v,w)(s)ds. \]
Using the linear estimates (\ref{107}) and (\ref{108}) we obtain
\begin{equation}\label{509}
\n{\Lambda_v w}{\X{\rho}{0}{\frac{1}{\rho}+}(\delta)} \le c \n{u_{\le}}{\widehat{L^{\rho}_{x}}}  + c \delta^{\frac{1}{\rho'}-} \n{N(v,w)}{\widehat{L^{\rho}_{xt}}(\delta)}
\end{equation}
with
\begin{eqnarray}\label{510}
\n{N(v,w)}{\widehat{L^{\rho}_{xt}}(\delta)} & \le & 2\n{|v|^2w}{\widehat{L^{\rho}_{xt}}(\delta)}+\n{v^2\overline{w}}{\widehat{L^{\rho}_{xt}}(\delta)} \\
& + & 2\n{v|w|^2}{\widehat{L^{\rho}_{xt}}(\delta)}+\n{\overline{v}|w|^2}{\widehat{L^{\rho}_{xt}}(\delta)}+\n{|w|^2w}{\widehat{L^{\rho}_{xt}}(\delta)}. \nonumber
\end{eqnarray}
Using part i) of Lemma \ref{l50} and Lemma \ref{l51} we get
\begin{equation}\label{511}
2\n{|v|^2w}{\widehat{L^{\rho}_{xt}}(\delta)}+\n{v^2\overline{w}}{\widehat{L^{\rho}_{xt}}(\delta)} \le c \delta^{\frac{1}{\rho}-\frac{1}{2}-} \q{v}{\XX{0}{\frac{1}{2}+}(\delta)}\n{w}{\X{\rho}{0}{\frac{1}{\rho}+}(\delta)},
\end{equation}
while Corollary \ref{k22} gives
\begin{equation}\label{512}
\n{|w|^2w}{\widehat{L^{\rho}_{xt}}(\delta)} \le c \|w\|^3_{\X{\rho}{0}{\frac{1}{\rho}+}(\delta)}.
\end{equation}
Estimates (\ref{511}) and (\ref{512}) are valid for $2 \ge \rho >1$. Although they are somehow intermediate, it is more complicated to discuss the 
quadratic-in-$w$-terms in (\ref{510}). We distinguish two cases:

\vspace{0,5cm}

{\bf {\underline{Case 1:}}} $r \ge \rho > \frac{4}{3}$. \\

In this case, part ii) of Lemma \ref{l50} combined with Lemma \ref{l51} gives 
\begin{equation}\label{513}
2\n{v|w|^2}{\widehat{L^{\rho}_{xt}}(\delta)}+\n{\overline{v}|w|^2}{\widehat{L^{\rho}_{xt}}(\delta)} \le c \delta^{(\frac{1}{4}-\frac{1}{2 \rho'})-}\n{v}{\XX{0}{\frac{1}{2}+}(\delta)}\q{w}{\X{\rho}{0}{\frac{1}{\rho}+}(\delta)}.
\end{equation}
Collecting the information from (\ref{509}) to (\ref{513}) we arrive at 
\[\n{\Lambda_v w}{\X{\rho}{0}{\frac{1}{\rho}+}(\delta)}  \le  c \n{u_{\le}}{\widehat{L^{\rho}_{x}}}  + c (\delta^{\frac{1}{2}-}\q{v}{\XX{0}{\frac{1}{2}+}(\delta)}+\delta^{\frac{1}{\rho'}-}\q{w}{\X{\rho}{0}{\frac{1}{\rho}+}(\delta)})\n{w}{\X{\rho}{0}{\frac{1}{\rho}+}(\delta)} \]
\begin{equation}\label{514}
 \le  c N^{\frac{r'}{\rho'}-1} + c (\delta^{\frac{1}{2}-}N^{r'-2} + \delta^{\frac{1}{\rho'}-}\q{w}{\X{\rho}{0}{\frac{1}{\rho}+}(\delta)})\n{w}{\X{\rho}{0}{\frac{1}{\rho}+}(\delta)}.
\end{equation}
Similarly we derive
\begin{eqnarray}\label{515}
&&\n{\Lambda_v w_1 -\Lambda_v w_2}{\X{\rho}{0}{\frac{1}{\rho}+}(\delta)}  \\
\le & c &(\delta^{\frac{1}{2}-}N^{r'-2}  + \delta^{\frac{1}{\rho'}-}(\q{w_1}{\X{\rho}{0}{\frac{1}{\rho}+}(\delta)}+\q{w_2}{\X{\rho}{0}{\frac{1}{\rho}+}(\delta)}))\n{w_1-w_2}{\X{\rho}{0}{\frac{1}{\rho}+}(\delta)}. \nonumber
\end{eqnarray}
Now choosing $R=2cN^{\frac{r'}{\rho'}-1}$ and $\delta$ such that $C\delta^{\frac{1}{2}-}N^{r'-2}=1$ (with a large constant $C$) we see that $\Lambda_v$
is a contraction of the closed ball of radius $R$ in $\X{\rho}{0}{\frac{1}{\rho}+}(\delta)$ into itself. By the contraction mapping principle we obtain
a unique solution $w \in \X{\rho}{0}{\frac{1}{\rho}+}(\delta)$ of (\ref{505}) with lifespan $\delta$ according to (\ref{503}). This solution satisfies (\ref{506}).

\vspace{0,5cm}

{\bf {\underline{Case 2:}}} $\frac{4}{3} \ge \rho > 1$. \\

Here we fix $\rho_0$ with
\begin{equation}\label{516}
\frac{3}{4}>\frac{1}{\rho_0}>\frac{3}{4}-\frac{1}{2\rho'}.
\end{equation}
Then, by the discussion concerning case 1, the estimates (\ref{514}) and (\ref{515}) hold with $\rho$ replaced by $\rho_0$. Using part iii) of
Lemma \ref{l50} we get
\begin{equation}\label{517}
2\n{v|w|^2}{\widehat{L^{\rho}_{xt}}(\delta)}+\n{\overline{v}|w|^2}{\widehat{L^{\rho}_{xt}}(\delta)} \le c \delta^{(\frac{5}{4}-\frac{3}{2 \rho'}-\frac{1}{\rho_0})-}\n{v}{\XX{0}{\frac{1}{2}+}(\delta)}\q{w}{\X{\rho_0}{0}{\frac{1}{\rho_0}+}(\delta)}
\end{equation}
instead of (\ref{513}). This gives, as substitute for (\ref{514}),
\begin{eqnarray*}
\n{\Lambda_v w}{\X{\rho}{0}{\frac{1}{\rho}+}(\delta)} & \le & c \n{u_{\le}}{\widehat{L^{\rho}_{x}}}  + c (\delta^{\frac{1}{2}-}\q{v}{\XX{0}{\frac{1}{2}+}(\delta)}+\delta^{\frac{1}{\rho'}-}\q{w}{\X{\rho}{0}{\frac{1}{\rho}+}(\delta)})\n{w}{\X{\rho}{0}{\frac{1}{\rho}+}(\delta)}  \\
&+&c \delta^{\alpha}\n{v}{\XX{0}{\frac{1}{2}+}(\delta)}\q{w}{\X{\rho_0}{0}{\frac{1}{\rho_0}+}(\delta)} 
\end{eqnarray*}
\begin{eqnarray}\label{518}
\hspace{1,5cm}&\le &    c N^{\frac{r'}{\rho'}-1} + c (\delta^{\frac{1}{2}-}N^{r'-2} + \delta^{\frac{1}{\rho'}-}\q{w}{\X{\rho}{0}{\frac{1}{\rho}+}(\delta)})\n{w}{\X{\rho}{0}{\frac{1}{\rho}+}(\delta)}   \\
\hspace{1,5cm}&+&c \delta^{\alpha}N^{\frac{r'}{2}-1}\q{w}{\X{\rho_0}{0}{\frac{1}{\rho_0}+}(\delta)}, \nonumber
\end{eqnarray}
where $\alpha = (\frac{5}{4}-\frac{1}{2 \rho'}-\frac{1}{\rho_0})-$. For the difference $\Lambda_v w_1 -\Lambda_v w_2$ we obtain
\begin{eqnarray}\label{519}
&&\n{\Lambda_v w_1 -\Lambda_v w_2}{\X{\rho}{0}{\frac{1}{\rho}+}(\delta)}  \\
\le & c &(\delta^{\frac{1}{2}-}N^{r'-2}  + \delta^{\frac{1}{\rho'}-}(\q{w_1}{\X{\rho}{0}{\frac{1}{\rho}+}(\delta)}+\q{w_2}{\X{\rho}{0}{\frac{1}{\rho}+}(\delta)}))\n{w_1-w_2}{\X{\rho}{0}{\frac{1}{\rho}+}(\delta)} \nonumber \\
+ & c &\delta^{\alpha}N^{\frac{r'}{2}-1}(\n{w_1}{\X{\rho_0}{0}{\frac{1}{\rho_0}+}(\delta)}+\n{w_2}{\X{\rho_0}{0}{\frac{1}{\rho_0}+}(\delta)})\n{w_1-w_2}{\X{\rho_0}{0}{\frac{1}{\rho_0}+}(\delta)}.\nonumber
\end{eqnarray}

At this point we introduce the complete metric space $(B,d)$, where
\[B=\{w \in \X{\rho}{0}{\frac{1}{\rho}+}(\delta) \cap \X{\rho_0}{0}{\frac{1}{\rho_0}+}(\delta): \n{w}{\X{\rho}{0}{\frac{1}{\rho}+}(\delta)} \le R,\n{w}{\X{\rho_0}{0}{\frac{1}{\rho_0}+}(\delta)} \le R_0\},\]
with $R=2cN^{\frac{r'}{\rho'}-1}$, $R_0=2cN^{\frac{r'}{\rho'_0}-1}$ and
\[d(w_1,w_2)=\n{w_1-w_2}{\X{\rho}{0}{\frac{1}{\rho}+}(\delta)}+\n{w_1-w_2}{\X{\rho_0}{0}{\frac{1}{\rho_0}+}(\delta)}.\]
Then, for $w, w_1, w_2 \in B$, we deduce from (\ref{514}), (\ref{515}) (with $\rho_0$ instead of $\rho$), (\ref{518}) and (\ref{519}) that
\begin{eqnarray*}
\n{\Lambda_v w}{\X{\rho_0}{0}{\frac{1}{\rho_0}+}(\delta)} & \le & \frac{R_0}{2} + c(\delta^{\frac{1}{2}-}N^{r'-2}+ \delta^{\frac{1}{\rho_0}-}R_0^2)R_0 ,\\
\n{\Lambda_v w}{\X{\rho}{0}{\frac{1}{\rho}+}(\delta)} & \le & \frac{R}{2} + c(\delta^{\frac{1}{2}-}N^{r'-2}+ \delta^{\frac{1}{\rho}-}R^2 + \delta^{\alpha}N^{\frac{r'}{2}-1}R_0^2R^{-1})R ,\\
d(\Lambda_v w_1,\Lambda_v w_2) & \le & c(\delta^{\frac{1}{2}-}N^{r'-2}+ \delta^{\frac{1}{\rho_0}-}R_0^2+ \delta^{\frac{1}{\rho}-}R^2 + \delta^{\alpha}N^{\frac{r'}{2}-1}R_0)d(w_1,w_2).
\end{eqnarray*}

Now for $\delta$ with $c\delta^{\frac{1}{2}-}N^{r'-2}= \frac{1}{4}$ we have
\begin{itemize}
\item[i)] $\delta^{\frac{1}{\rho'_0}-}R_0^2 \sim N^{\frac{4}{\rho'_0}-2}$,
\item[ii)] $\delta^{\frac{1}{\rho'}-}R^2 \sim N^{\frac{4}{\rho'}-2}$ and
\item[iii)] $\delta^{\alpha}N^{\frac{r'}{2}-1}R_0 \le \delta^{\alpha}N^{\frac{r'}{2}-1}R^2_0R^{-1}\sim N^{3-\frac{2}{\rho'}-\frac{4}{\rho_0}}$.
\end{itemize}
All the exponents in i) - iii) are negative (concerning iii) cf. (\ref{516})), so that for $N$ sufficiently large the mapping $\Lambda_v$ becomes a contraction 
of $(B,d)$ into itself. Hence, we get a solution $w \in \X{\rho}{0}{\frac{1}{\rho}+}(\delta)\cap\X{\rho_0}{0}{\frac{1}{\rho_0}+}(\delta)$ of (\ref{505}),
being unique in $\X{\rho_0}{0}{\frac{1}{\rho_0}+}(\delta)$ and satisfying (\ref{506}). Again, the lifespan $\delta$ of this solution is given by (\ref{503}).

\subsection{Growth bounds for the $L^2$-norm of the regular part} \hfill \\

In order to extend the local solution $u$ of (\ref{101}) to a given time interval $[0,T]$, where $T$ is arbitrarily large, we shall glue together
local solutions of lifespan $\delta \sim N^{4-2r'-}$ until $T$ is reached. After a first step, at time $\delta$, the following Cauchy problems are
considered:
\begin{eqnarray*}
iv'_t + v'_{xx}+ |v'|^2v'=0 &;& v'(0)=v(\delta)+y(\delta)\\
iw'_t + w'_{xx}+ N(v',w')=0 &;& w'(0)=e^{i\delta \partial^2}u_{\le}.
\end{eqnarray*}

Here $v$ and $w$ are the local solutions of (\ref{502}) and (\ref{505}), respectively, living on $[0,\delta]$; $N(v,w)$ is as introduced below
(\ref{505}) and
\[y(t)= i \int_0^te^{i(t-s) \partial^2}N(v,w)(s)ds.\]
In order to reapply the local results concerning (\ref{502}) and (\ref{505}) at time $\delta$, we first observe that $e^{i\delta \partial^2}u_{\le}$
fulfills (\ref{500}). Then we have to make sure, that $y(\delta)$ belongs to $L^2_x$ and obeys (\ref{501}). Moreover, when repeating the argument
$\frac{T}{\delta}$ times until $T$ is reached, the total increment of the $L^2_x$ - norm of the regular (i. e. the $v$-) part of the solution $u$ must not
exceed its size at the beginning given by (\ref{501}). So we have to estimate
\[\sup_{0 \le t \le \delta} \n{y(t)}{L^2_x} \le c \n{y}{\XX{0}{\frac{1}{2}+ (\delta)}}\le c \n{N(v,w)}{\XX{0}{-\frac{1}{2}+ }(\delta)}:\]

\begin{lemma}\label{l52}
For any $\rho \in (1,2]$ and for $\frac{1}{\rho_0}=\frac{1}{4}+\frac{1}{2 \rho}$ the following estimate holds true:
\begin{eqnarray*}
\n{N(v,w)}{\XX{0}{-\frac{1}{2}+} (\delta)} \hspace{9cm} \\
 \le  c ( \delta ^{\frac{1}{4}+\frac{1}{2 \rho '}-} \q{v}{\XX{0}{\frac{1}{2}+} (\delta)} + \n{v}{\XX{0}{\frac{1}{2}+} (\delta)}\n{w}{\X{\rho}{0}{\frac{1}{\rho}+} (\delta)} +\q{w}{\X{\rho_0}{0}{\frac{1}{\rho_0}+} (\delta)})\n{w}{\X{\rho}{0}{\frac{1}{\rho}+} (\delta)}.
\end{eqnarray*}
\end{lemma}

Proof: In Lemma \ref{l21} we choose $q=2$, $\frac{1}{p'}=\frac{1}{4}+\e$, $\frac{1}{r'_0}=\frac{3}{4}-2\e$, $\frac{1}{r'_1}=\frac{1}{4}$ and 
$\frac{1}{r'_2}=4\e$. Using the symmetry between the first two factors and bilinear interpolation (with $\theta = \frac{1}{2}$) we see that 
for $u$, $v$ and $w$ as in Lemma \ref{l21}
\[\n{uvw}{L^2_x(\widehat{L^p_t})} \le c \n{u_0}{\widehat{L^{\rho_0}_x}}\n{v_0}{\widehat{L^{\rho_{0}}_x}}\n{w_0}{\widehat{L^{\rho_{1}}_x}},\]
provided $\frac{1}{p}=\frac{3}{4}-\e$, $\frac{1}{\rho_0}=\frac{1}{2}+\e$ and $\frac{1}{\rho_1}=1-4\e$. On the other hand, by H\"older and 
(\ref{102})
\[\n{uvw}{L^2_{xt}}\le c \n{u_0}{\widehat{L^{q_0}_x}}\n{v_0}{\widehat{L^{q_{0}}_x}}\n{w_0}{\widehat{L^{q_{1}}_x}},\]
where $\frac{1}{q_0}=\frac{3}{8}+\frac{\e}{2}$, $\frac{1}{q_1}=\frac{3}{4}-\e$. Interpolating again with $\theta$ chosen such that $\e \theta = 2\e- \frac{\theta}{4}$
leads to
\begin{equation}\label{520}
\n{uvw}{L^2_x(\widehat{L^r_t})} \le c \n{u_0}{L^2_x}\n{v_0}{L^2_x}\n{w_0}{\widehat{L^{\rho}_x}},
\end{equation}
whenever $2 \ge r > \frac{4}{3}$ and $\frac{2}{r}=\frac{1}{2}+\frac{1}{\rho}$. The corresponding $\x$ - estimate reads
\begin{equation}\label{521}
\n{fg\overline{h}}{L^2_x(\widehat{L^r_t})} \le c \n{f}{\XX{0}{\frac{1}{2}+}}\n{g}{\XX{0}{\frac{1}{2}+}}\n{h}{\X{\rho}{0}{\frac{1}{\rho}+}}.
\end{equation}
Combining this with Lemma \ref{l51} we obtain for the $v^2\overline{w}$ - term in $N(v,w)$:
\begin{eqnarray}\label{522}
\n{v^2\overline{w}}{\XX{0}{-\frac{1}{2}+}(\delta)} & \le & c \delta ^{b + \frac{1}{2}-} \n{v^2\overline{w}}{\XX{0}{b}(\delta)}\hspace{1cm}(0>b>-\frac{1}{2}) \nonumber\\
& \le & c \delta ^{ \frac{1}{r'}-} \n{v^2\overline{w}}{L^2_x(\widehat{L^r_t})(\delta)}\hspace{1cm}(\frac{1}{r}=(\frac{1}{2}-b)-)\nonumber\\
& \le & c \delta ^{\frac{1}{4}+\frac{1}{2 \rho '}-} \q{v}{\XX{0}{\frac{1}{2}+} (\delta)}\n{w}{\X{\rho}{0}{\frac{1}{\rho}+} (\delta)}.
\end{eqnarray}
To treat the $v|w|^2$ - term in $N(v,w)$ we use Corollary \ref{k22} directly:
\begin{equation}\label{523}
\n{v|w|^2}{\XX{0}{-\frac{1}{2}+}(\delta)} \le c \n{v|w|^2}{L^2_x(\widehat{L^{\rho}_t})(\delta)} \le c \n{v}{\XX{0}{\frac{1}{2}+} (\delta)}\q{w}{\X{\rho}{0}{\frac{1}{\rho}+} (\delta)}.
\end{equation}
From Corollary \ref{k22} we conclude further - by symmetry between the first two factors and bilinear interpolation - that for $\frac{1}{\rho_0}=\frac{1}{4}+\frac{1}{2 \rho}$
\[\n{fg\overline{h}}{L^2_x(\widehat{L^{\rho}_t})} \le c \n{f}{\X{\rho_0}{0}{\frac{1}{\rho_0}+}}\n{g}{\X{\rho_0}{0}{\frac{1}{\rho_0}+}}\n{h}{\X{\rho}{0}{\frac{1}{\rho}+}}.\]
This gives
\begin{equation}\label{524}
\n{|w|^2w}{\XX{0}{-\frac{1}{2}+}(\delta)} \le c \n{|w|^2w}{L^2_x(\widehat{L^{\rho}_t})(\delta)} \le c \q{w}{\X{\rho_0}{0}{\frac{1}{\rho_0}+} (\delta)}\n{w}{\X{\rho}{0}{\frac{1}{\rho}+} (\delta)}.
\end{equation}
It remains to consider the two terms containing a factor $\overline{v}$. In order to treat the $w|v|^2$ - contribution in $N(v,w)$, we go back to
(\ref{521}), that is - after replacing $\rho$ by $\rho_1$ -
\[\n{fg\overline{h}}{L^2_x(\widehat{L^r_t})} \le c \n{f}{\XX{0}{\frac{1}{2}+}}\n{g}{\XX{0}{\frac{1}{2}+}}\n{h}{\X{\rho_1}{0}{\frac{1}{\rho_1}+}},\]
$\rho_1 > 1$, $\frac{1}{r}=\frac{1}{4}+\frac{1}{2 \rho_1}$, telling us that
\[M_{g\overline{h}}: \XX{0}{\frac{1}{2}+} \longrightarrow L^2_x(\widehat{L^r_t});\hspace{1cm}f \mapsto fg\overline{h}\]
is continuous with operator norm bounded by $c\n{g}{\XX{0}{\frac{1}{2}+}}\n{h}{\X{\rho_1}{0}{\frac{1}{\rho_1}+}}$. But then the adjoint operator
\[M_{\overline{g}h}: L^2_x(\widehat{L^{r'}_t}) \longrightarrow \XX{0}{-\frac{1}{2}-};\hspace{1cm}f \mapsto f\overline{g}h\]
is also bounded with the same norm, which gives us the estimate (after exchanging $g$ and $h$)
\[\n{fg\overline{h}}{\XX{0}{-\frac{1}{2}-}} \le c \n{f}{L^2_x(\widehat{L^{r'}_t})}\n{g}{\X{\rho_1}{0}{\frac{1}{\rho_1}+}}\n{h}{\XX{0}{\frac{1}{2}+}} \le c \n{f}{\XX{0}{b_1}}\n{g}{\X{\rho_1}{0}{\frac{1}{\rho_1}+}}\n{h}{\XX{0}{\frac{1}{2}+}},\]
whenever $b_1 > \frac{1}{2 \rho_1} - \frac{1}{4}$. Interpolation with the $L^6_{xt}$ - Strichartz estimate in the form
\[\n{fg\overline{h}}{L^2_{xt}} \le c \n{f}{\XX{0}{\frac{1}{2}+}}\n{g}{\XX{0}{\frac{1}{2}+}}\n{h}{\XX{0}{\frac{1}{2}+}}\]
gives
\[\n{fg\overline{h}}{\XX{0}{-\frac{1}{2}+}} \le c \n{f}{\XX{0}{b}}\n{g}{\X{\rho}{0}{\frac{1}{\rho}+}}\n{h}{\XX{0}{\frac{1}{2}+}},\]
where now $b > \frac{1}{2 \rho} - \frac{1}{4}$ is necessary. Using Lemma \ref{l51} again we obtain
\begin{equation}\label{525}
\n{w|v|^2}{\XX{0}{-\frac{1}{2}+}(\delta)}  \le c \delta ^{\frac{1}{4}+\frac{1}{2 \rho '}-} \q{v}{\XX{0}{\frac{1}{2}+} (\delta)}\n{w}{\X{\rho}{0}{\frac{1}{\rho}+} (\delta)}
\end{equation}
in close analogy to (\ref{522}).
In order to prove
\begin{equation}\label{526}
\n{w^2\overline{v}}{\XX{0}{-\frac{1}{2}+}(\delta)}  \le c \n{v}{\XX{0}{\frac{1}{2}+} (\delta)}\q{w}{\X{\rho}{0}{\frac{1}{\rho}+} (\delta)}.
\end{equation}
we start from
\begin{equation}\label{527}
\n{I^{\frac{1}{2}}(f \overline{g})}{\widehat{L^{\rho_1}_x}(L^2_t)} \le c \n{f}{\XX{0}{\frac{1}{2}+}}\n{g}{\X{\rho_1}{0}{\frac{1}{\rho_1}+}},
\end{equation}
where $2 \ge \rho_1 >1$ (Corollary \ref{k20}). Interpolation with
\[\n{f \overline{g}}{L^2_{xt}} \le \n{f }{L^3_{xt}}\n{g}{L^6_{xt}} \le c \n{f }{\XX{0}{\frac{1}{4}+}}\n{g}{\XX{0}{\frac{1}{2}+}}\]
(Strichartz, for the first factor interpolated with $L^2_{xt} = \XX{0}{0}$) yields
\[\n{I^s(f \overline{g})}{\widehat{L^{\rho}_x}(L^2_t)} \le c \n{f}{\XX{0}{b}}\n{g}{\X{\rho}{0}{\frac{1}{\rho}+}},\]
provided $2 \ge \rho >1$, $\frac{1}{2} \ge s > \frac{1}{\rho} - \frac{1}{2}$, $b > \frac{1}{4} + \frac{s}{2}$. Dualizing we obtain
\begin{equation}\label{528}
\n{gh}{\XX{0}{-b}} \le c \n{g}{\X{\rho}{0}{\frac{1}{\rho}+}}\n{I^{-s}h}{\widehat{L^{\rho'}_x}(L^2_t)}.
\end{equation}
Now for $f, g \in \X{\rho}{0}{\frac{1}{\rho}+}$, $h \in \XX{0}{\frac{1}{2}+}$ and $P= \F^{-1} \chi_{|\xi|>1}\F$ we get from (\ref{528}) and 
(\ref{527})
\begin{eqnarray}\label{529}
\n{f P(g \overline{h})}{\XX{0}{-\frac{1}{2}+}} & \le & c \n{f}{\X{\rho}{0}{\frac{1}{\rho}+}}\n{I^{-\frac{1}{2}+}P(g \overline{h})}{\widehat{L^{\rho'}_x}(L^2_t)} \nonumber \\
& \le & c \n{f}{\X{\rho}{0}{\frac{1}{\rho}+}}\n{I^{\frac{1}{2}}(g \overline{h})}{\widehat{L^{\rho}_x}(L^2_t)} \\
& \le & c \n{f}{\X{\rho}{0}{\frac{1}{\rho}+}}\n{g}{\X{\rho}{0}{\frac{1}{\rho}+}}\n{h}{\XX{0}{\frac{1}{2}+}}. \nonumber
\end{eqnarray}
By symmetry between $f$ and $g$ we have the same bound for $\n{g P(f \overline{h})}{\XX{0}{-\frac{1}{2}+}}$. To establish (\ref{526}) it remains to 
remove the projector $P$ from (\ref{529}). For that purpose we choose functions $\psi_1$, $\psi_2$, $\psi_3$ with 
\[\n{\psi_1}{L^{\rho'}_{\xi,\tau}}= \n{f}{\X{\rho}{0}{\frac{1}{\rho}+}},\,\,\,\n{\psi_2}{L^{\rho'}_{\xi,\tau}}= \n{g}{\X{\rho}{0}{\frac{1}{\rho}+}},\,\,\,\n{\psi_3}{L^{2}_{\xi,\tau}}= \n{h}{\XX{0}{\frac{1}{2}+}}.\]
Then, assuming $\hat{f}$, $\hat{g}$ and $\hat{h}$ to be nonegative,
\begin{eqnarray*}
\F(f g \overline{h})(\xi,\tau)  \le \F(f P(g \overline{h}))(\xi,\tau)  +\F( g P(f \overline{h}))(\xi,\tau) \hspace{4cm}\\
+ c \int d\nu \psi_1(\xi_1,\tau_1)\langle \sigma_1 \rangle ^{-\frac{1}{\rho}-}\chi(\xi-\xi_1)\psi_2(\xi_2,\tau_2)\langle \sigma_2 \rangle ^{-\frac{1}{\rho}-}\chi(\xi-\xi_2)\psi_3(\xi_3,\tau_3)\langle \sigma_3 \rangle ^{-\frac{1}{2}-},
\end{eqnarray*}
where $d\nu= d\xi_1d\tau_1d\xi_2d\tau_2$, $\sum_{i=1}^3(\tau_i , \xi_i) =(\tau, \xi)$, $\sigma_{1,2} = \tau_{1,2} + \xi_{1,2} ^2$, $\sigma_3 = \tau_3 - \xi_3 ^2$, $\chi=\chi_{[-1,1]}$.
Writing $I(\xi,\tau)$ for the above integral, it is sufficient to show that
\begin{equation}\label{530}
\n{I}{L^{2}_{\xi,\tau}} \le c \n{\psi_1}{L^{\rho'}_{\xi,\tau}}\n{\psi_2}{L^{\rho'}_{\xi,\tau}}\n{\psi_3}{L^{2}_{\xi,\tau}}.
\end{equation}
Throwing away the $\sigma_3$ - factor and using H\"older's inequality we have
\begin{eqnarray*}
&&|I(\xi,\tau)| \\
& \le & \n{\psi_1}{L^{\rho'}_{\xi,\tau}} \n{\psi_2}{L^{\rho'}_{\xi,\tau}}\left(\int d\nu \chi(\xi-\xi_1)\chi(\xi-\xi_2)\langle \sigma_1 \rangle ^{-1-}\langle \sigma_2 \rangle ^{-1-} |\psi_3(\xi_3,\tau_3)|^{\rho}\right)^{\frac{1}{\rho}}\\
& \le & c \n{\psi_1}{L^{\rho'}_{\xi,\tau}} \n{\psi_2}{L^{\rho'}_{\xi,\tau}}\left(\int d\nu \chi(\xi-\xi_1)\chi(\xi-\xi_2)\langle \sigma_1 \rangle ^{-1-}\langle \sigma_2 \rangle ^{-1-} |\psi_3(\xi_3,\tau_3)|^{2}\right)^{\frac{1}{2}}.
\end{eqnarray*}
Hence
\begin{eqnarray*}
&&\n{I}{L^{2}_{\xi,\tau}} \\
&\le & c\n{\psi_1}{L^{\rho'}_{\xi,\tau}} \n{\psi_2}{L^{\rho'}_{\xi,\tau}}\n{(\int d\nu \chi(\xi-\xi_1)\chi(\xi-\xi_2)\langle \sigma_1 \rangle ^{-1-}\langle \sigma_2 \rangle ^{-1-} |\psi_3(\xi_3,\tau_3)|^{2})^{\frac{1}{2}}}{L^{2}_{\xi,\tau}}
\end{eqnarray*}
where the square of the last factor is equal to
\[\int d\xi d\tau d \nu \chi(\xi-\xi_1)\chi(\xi-\xi_2)\langle \sigma_1 \rangle ^{-1-}\langle \sigma_2 \rangle ^{-1-} |\psi_3(\xi_3,\tau_3)|^{2} \le c \q{\psi_3}{L^{2}_{\xi,\tau}}\]
as desired. Now (\ref{530}) and thus (\ref{526}) are shown. Finally, collecting the information from (\ref{522}) - (\ref{526}), we obtain the claimed estimate.
$\hfill \Box$

\vspace{0,3cm}

Now taking into account that $\delta \sim N^{(4-2r')-}$, $\n{v}{\XX{0}{\frac{1}{2}+} (\delta)} \le c N^{\frac{r'}{2}-1}$ and, for any $\rho  \in (1,r]$,
$\n{w}{\X{\rho}{0}{\frac{1}{\rho}+} (\delta)} \le c N^{\frac{r'}{\rho'}-1}$ (see (\ref{503}), (\ref{504}) and (\ref{506})), we see that for $2 \ge \rho > 1$,
$\frac{1}{\rho_0}=\frac{1}{4}+\frac{1}{2 \rho}$
\begin{eqnarray*}
\sup_{0 \le t \le \delta} \n{y(t)}{L^2_x} & \le & c (N^{(\frac{r'}{2}-1)+} + N^{\frac{r'}{2}+\frac{r'}{\rho'}-2}+N^{\frac{2r'}{\rho'_0}-2})N^{\frac{r'}{\rho'}-1}\\
& \le & c N^{(\frac{r'}{2}+\frac{r'}{\rho'}-2)+},
\end{eqnarray*}
which is the bound for the growth of the $L^2_x$ - norm of the $v$ - part in each step. Setting $N^{\e}=T$ (for some $\e = \e (r) > 0$ to be specified below)
the number of iterations necessary to reach $T$ becomes
\[\frac{T}{\delta} \sim N^{(2r'-4+\e )+},\]
giving a total increment of mass of about
\[N^{(\frac{5r'}{2}+\frac{r'}{\rho'}-6+\e) +} \le N^{\frac{r'}{2}-1},\]
the latter, provided
\begin{equation}\label{531}
2r'-5+\e+\frac{r'}{\rho'} < 0.
\end{equation}
Choosing $\e = 5 - 2r'-$ and $\rho$ close to $1$ so that (\ref{531}) is fulfilled, we have
\[N= T^{\frac{1}{\e}}=T^{\frac{1+}{5-2r'}}\]
and the iteration process described above yields a solution $u$ of (\ref{101}) defined on the whole interval $[0,T]$ and satisfying
\[\n{u(T)-e^{iT\partial ^2}u_0}{L^2_x} \le c N^{\frac{r'}{2}-1} = c T^{\frac{r'-2}{10-4r'}+}.\]
This concludes the proof of Theorem \ref{t4}.

\end{document}